\documentclass[a4paper,11pt,reqno]{amsart}

\usepackage{latexsym,amssymb,dsfont,amsmath,amsthm}

\chardef\bslash=`\\ 

\numberwithin{equation}{section}

\newtheorem{theorem}{Theorem}[section]
\newtheorem{corollary}[theorem]{Corollary}
\newtheorem{lemma}[theorem]{Lemma}
\newtheorem{proposition}[theorem]{Proposition}

\theoremstyle{remark}
\newtheorem{remark}[theorem]{Remark}
\newtheorem{example}[theorem]{Example}
\theoremstyle{definition}
\newtheorem{definition}[theorem]{Definition}

\newcommand{\thmref}[1]{Theorem~\ref{#1}}
\newcommand{\secref}[1]{Section~\ref{#1}}
\newcommand{\proref}[1]{Proposition~\ref{#1}}
\newcommand{\lemref}[1]{Lemma~\ref{#1}}

\newcommand\bp{\begin{proof}}
\newcommand\ep{\end{proof}}

\newcommand\3[1]{{\mathds #1}}

\newcommand{\N}{\mathbb N}
\newcommand{\Z}{\mathbb Z}
\newcommand{\HH}{\mathbb H}
\newcommand{\Q}{\mathbb Q}
\newcommand{\R}{\mathbb R}
\newcommand{\C}{\mathbb C}

\newcommand\af{\mathbb{A}_f}

\newcommand{\hh}{\mathcal H}
\newcommand{\rr}{\hat{\Z}}
\newcommand{\primes}{\mathcal  P}

\newcommand\bpmatrix{\begin{pmatrix}}
\newcommand\epmatrix{\end{pmatrix}}

\newcommand{\matr}[2]{\left(\begin{matrix}1&#1 \\
0&#2\end{matrix}\right)}
\newcommand{\diag}[2]{\left(\begin{matrix}#1&0\\
0&#2\end{matrix}\right)}

\newcommand\mz{{\operatorname{Mat}_2}(\Z)}
\newcommand\mzp{{\operatorname{Mat}^+_2}(\Z)}
\newcommand\mq{{\operatorname{Mat}_2}(\Q)}
\newcommand\mr{{\operatorname{Mat}_2}(\rr)}
\newcommand\ma{{\operatorname{Mat}_2}(\af)}
\newcommand\mtwo{\operatorname{Mat}_2}
\newcommand\mn{{\operatorname{Mat}_n}}

\newcommand\HA{\operatorname{HAut}}

\newcommand\glq{{\operatorname{GL}^+_2}(\Q)}
\newcommand\gla{{\operatorname{GL}_2}(\af)}
\newcommand\glap{{\operatorname{GL}^+_2}(\af)}
\newcommand\slz{{\operatorname{SL}_2}(\Z)}
\newcommand\slr{{\operatorname{SL}_2}(\rr)}
\newcommand\glr{{\operatorname{GL}_2}(\rr)}
\newcommand\gl{{\operatorname{GL}_2}}
\newcommand\glp{{\operatorname{GL}_2^+}}
\newcommand\sltwo{{\operatorname{SL}_2}}
\newcommand\gln{{\operatorname{GL}_n}}
\newcommand\glf{{\operatorname{GL}_4}}
\newcommand\glo{{\operatorname{GL}_1}}
\newcommand\G{\gl}

\newcommand{\Ad}{\operatorname{Ad}}
\newcommand{\Aut}{\operatorname{Aut}}
\newcommand{\id}{\operatorname{id}}
\newcommand{\Ker}{\operatorname{Ker}}
\newcommand{\lspa}{\operatorname{span}}
\newcommand{\Tr}{\operatorname{Tr}}
\newcommand{\supp}{\operatorname{supp}}

\newcommand{\inv}{^{-1}}

\newcommand{\hecke}[2]{\mathcal \hh({#1}, {#2})}
\newcommand{\redheck}[2]{C^*_r({#1}, {#2})}
\newcommand{\red}[1]{C^*_r(#1)}

\newcommand{\fib}[3]{#2\backslash #1\times_{#2}#3}
\newcommand{\fibb}[3]{#2\backslash #1\boxtimes_{#2}#3}

\newcommand\enu[1]{\smallskip\newline\makebox[5mm][l]{\rm(#1)}}

\begin{document}

\title[Hecke algebras]{Hecke algebras of semidirect products and
the finite part of the Connes-Marcolli C$^*$-algebra}

\thanks{
Part of this work was carried out through several visits of the
first author to the Department of Mathematics at the University of
Oslo. He would like to thank the department for their hospitality,
and the SUPREMA project for the support. }

\author[M. Laca]{Marcelo Laca$^1$}
\address{Department of Mathematics and Statistics, University of
Victoria, PO Box 3045, Victoria, British Columbia, V8W 3P4, Canada.}
\email{laca@math.uvic.ca}
\thanks{$^1$) Supported by the Natural Sciences and Engineering
Research Council of Canada. }

\author[N. S. Larsen]{Nadia S. Larsen$^2$}
\address{Department of Mathematics, University of Oslo,
P.O. Box 1053 Blindern, N-0316 Oslo, Norway.}
\email{nadiasl@math.uio.no}
\thanks{$^2$) Supported by the  Research Council of Norway. }

\author[S. Neshveyev]{Sergey Neshveyev$^2$}
\address{Department of Mathematics, University of Oslo,
P.O. Box 1053 Blindern, N-0316 Oslo, Norway.}
\email{sergeyn@math.uio.no}

\date{October 16, 2006; minor corrections June 20, 2007.}

\begin{abstract}
We study a C$^*$-dynamical system arising from the ring inclusion of
the $2\times 2$ integer matrices in the rational ones. The
orientation preserving affine groups of these rings form a Hecke
pair that is closely related to a recent construction of Connes and
Marcolli; our dynamical system consists of the associated reduced
Hecke C$^*$-algebra endowed with a canonical dynamics defined in
terms of the determinant function. We show that the  Schlichting
completion also consists of affine groups of matrices, over the
finite adeles, and we obtain results about the structure and induced
representations of the Hecke C$^*$-algebra. In a somewhat unexpected
parallel with the one dimensional case studied by Bost and Connes,
there is a group of symmetries given by an action of the finite
integral ideles, and the corresponding fixed point algebra
decomposes as a tensor product over the primes.  This decomposition
allows us to obtain a complete description of a natural class of
equilibrium states which conjecturally includes all
KMS$_\beta$-states for $\beta\ne0,1$.
\end{abstract}

\maketitle

\bigskip

\section*{Introduction}

Let $\HH$ denote the upper halfplane and let $\af$ be the ring of
finite adeles. The group $\glq$ of  $2\times 2$ rational matrices
with positive determinant acts on $\HH\times \ma$, by M\"obius
transformations on $\HH$ and by left multiplication on $\ma$.
Roughly speaking, the C$^*$-algebra underlying the $\G$-system of
Connes and Marcolli~\cite{con-mar} can be constructed by effecting
two modifications on the corresponding transformation groupoid
$\glq\times(\HH\times\ma)$; the first one is to cut down from~$\af$
to the compact open subring of  finite integral adeles
$\rr=\prod_p\Z_p$  yielding a kind of semigroup crossed product by
$\mzp$, and the second one is to eliminate the degeneracy due to the
$\Gamma = \slz$ symmetries by factoring out the action of
$\Gamma\times\Gamma$ given by $(\gamma_1, \gamma_2) (g,x) =
(\gamma_1 g \gamma_2\inv, \gamma_2 x)$. These modifications destroy
the initial groupoid and semigroup crossed product structures, but
the resulting C$^*$-algebra
$$
\red{\fibb{\glq}{\Gamma}{(\HH\times\mr)}},
$$
for which we use the notation of \cite{LLN}, retains enough of the
original transformation groupoid flavor that it is possible to use
slightly modified crossed product techniques in its study. The
convolution formula that gives the product on the Connes-Marcolli
C$^*$-algebra is based on the convolution formula for the classical
Hecke algebra. This connection was pursued early on by
Tzanev~\cite{tza2}, who pointed out that the Connes-Marcolli
C$^*$-algebra could also be described as $\red{\fib{P}{P_0}{\HH}}$,
where
$$
(P,P_0) = \left(\matr{\mq}{\glq}, \matr{\mz}{\slz} \right)
$$
is a Hecke pair and the action of $P$ on $\HH$ is defined by
M\"obius transformations through the obvious homomorphism $P\to\glq$
(technically,  the action of $P_0$ is not proper, but the
construction of $\red{\fib{P}{P_0}{\HH}}$ still makes sense
by~\cite[Remark~1.4]{LLN}). Thus, the Connes-Marcolli C$^*$-algebra
can be thought of as a new type of crossed product: of the algebra
$C_0(\Gamma\backslash\HH)$ by the Hecke pair $(P,P_0)$, and, in
particular, the reduced Hecke C$^*$-algebra $\redheck{P}{P_0}$ is
contained in the multiplier algebra of the Connes-Marcolli
C$^*$-algebra, see \cite[Lemma~1.3]{LLN}. Because of this, abusing
slightly the terminology, we will refer to $\redheck{P}{P_0}$ as the
{\em finite part of the Connes-Marcolli C$^*$-algebra}, and we point
out  that this finite part corresponds to the quotient of the
determinant part of the $\G$-system,  cf.  \cite[Section
1.7]{con-mar}, by the above action of $\Gamma \times \Gamma$.

The goal of the present work is to study the structure of
$\redheck{P}{P_0}$ and the phase transition of the corresponding
C$^*$-dynamical system. We were initially motivated by our belief
that it should be possible to exploit the crossed product structure
observed by Tzanev in order to study the phase transition of the
Connes-Marcolli system, and that in order to do this one would have
to understand first the structure and the phase transition of the
finite part. We were also motivated by the observation that the
Hecke pair $(P,P_0)$ consists of the orientation preserving affine
transformations of the rings of $2\times 2$ matrices over the
rationals and over the integers, and hence the associated
C$^*$-dynamical system is a very natural (albeit somewhat na\"{i}ve)
higher dimensional version of the one studied by Bost and
Connes~\cite{bos-con}, which certainly deserves consideration. In
addition $(P,P_0)$ is a very interesting example of a Hecke inclusion
of semidirect product groups, a class that has received considerable
attention in recent years, see e.g. \cite{bre,LaLa,larr,klq}.

As it turned out,  we were able to study the Connes-Marcolli phase
transition and to prove the uniqueness of the KMS$_\beta$-states for
$\beta$ in the critical interval by a more direct method that does
not require consideration of $\redheck{P}{P_0}$, although it does
rely on it for insight, see \cite{LLN}. Interestingly enough, the
phase transition of the finite part of the Connes-Marcolli system
seems to be a more resilient problem than for their full
$\G$-system. The main reason for this is that the freeness resulting
from the `infinite part', that is to say, freeness of $\glq$ acting
on $\HH\times(\ma\setminus\{0\})$, is a crucial ingredient in
reducing KMS-states of the Connes-Marcolli $\G$-system to measures
on $\HH\times \ma$. Because of this our classification of the
KMS-states of the finite part relies on an extra hypothesis of
regularity which allows us to use techniques similar to those of
\cite{bos-con,minautex,nes,con-mar,LLN}. This regularity property
seems natural and we believe it to be automatic, but we have not
been able to prove this.

A brief summary of the contents of each section follows. In
\secref{section_gen_P_P0} we study Hecke pairs of semidirect
products modeled on our main example $(P,P_0)$, but general enough
to be of independent interest. The results include necessary and
sufficient conditions for an inclusion of semidirect products to be
a Hecke pair, structural results that highlight the role of the
Hecke algebras of the factors in that of the semidirect product, and
the discussion of a (nonselfadjoint) Hecke algebra associated to the
Hecke inclusion of a group in a semigroup determined by the
inclusion of semidirect products. In \secref{completions} we study
completions of a given Hecke pair of semidirect products to
topological Hecke pairs. We show that a completion of a semidirect
product $V\rtimes G$ can be chosen to be a semidirect product
itself, which can be computed in terms of completions of $V$ and
$G$. We then use this fact to construct induced representations of
the Hecke algebra from representations of the group $V$, and we use
this in \thmref{pi_faithful_our_pair} to obtain certain faithful
representations of $\redheck{P}{P_0}$ indexed by $\glr$ in the case
of our main motivating example. In \secref{symmetries} we show that
the finite part of the $\G$-system carries an action of the group of finite integral ideles 
$\rr^*=\prod_p\Z_p^*$ as
symmetries, whose fixed point algebra decomposes as a tensor product
over the primes, \thmref{Symm}. Therefore the situation in this
finite part is rather surprisingly quite similar to the one
dimensional case of Bost and Connes~\cite{bos-con}. We also
characterize two important subalgebras  of this fixed point algebra
in terms of integer lattices. The results on representations,
completions and symmetries come together in \secref{SectionKMS},
where we obtain our classification result, \thmref{main_theorem}, in
which we show that the phase transition of regular KMS-states for
$\beta > 2$  is indexed by the symmetry group $\rr^*$, and that for
$\beta \in (1,2]$ there is a unique regular KMS$_\beta$-state.

\bigskip

\section{Hecke pairs from semidirect products}
\label{section_gen_P_P0}

Motivated by the inclusion of the orientation preserving affine
group of $\mz$ in that of $\mq$, we are interested in the following
general situation: let $G$ be a group acting by automorphisms of
another group $V$, and let $\Gamma$ be a subgroup of $G$ leaving a
subgroup $V_0$ of $V$ invariant,  so the semidirect product
$V_0\rtimes \Gamma$ can be viewed as a subgroup of $V\rtimes G$. We
aim to study the Hecke algebra of the inclusion $V_0\rtimes\Gamma
\subset V\rtimes G$, and start by recalling basic definitions.

\subsection{Hecke pairs and their $*$-algebras}
\label{subsection_hecke_alg} As customary, by a {\em Hecke pair}
$(N, N_0)$ we mean a group~$N$ with a subgroup~$N_0$ such that
$$
L_{N_0}(x): =[N_0:N_0\cap x N_0x^{-1}]
$$
is finite for all $x$ in $N$.  More generally, when $X$ is  a right
$N_0$-invariant subset of~$N$  we denote by~$L_{N_0}(X)$ the number
of left $N_0$-cosets in $X$. Then $L_{N_0}(x)=L_{N_0}(N_0xN_0)$, so
the quantity $L_{N_0}(x)$ is the number of left cosets in the double
coset $N_0 xN_0$, and $R_{N_0}(x):= L_{N_0}(x^{-1})$ is the number
of right cosets. We shall often drop the subindex from the notation
$L_{N_0}(x)$ when there is no risk of confusion.

The formula $\Delta_{N_0}(x):=L_{N_0}(x)/L_{N_0}(x^{-1})$ defines a
homomorphism $N\to \Q^+$, see for instance~\cite[Proposition
I.3.6]{kri}. We refer to it as the modular function of $(N,N_0)$.

The \emph{Hecke algebra} $\hh(N,N_0)$ of a Hecke pair $(N,N_0)$
consists of the vector space of  complex valued $N_0$-biinvariant
functions supported on finitely many double cosets,
$$
\{f\colon N\to \C\mid f(n_1xn_2)=f(x),\forall\ n_1,n_2\in N_0,\
\supp(f)\text{ finite in } N_0\backslash N/N_0\},
$$
endowed with the product
\begin{equation}
(f_1\ast f_2)(x)=\sum_{y\in N_0\backslash N}f_1(xy^{-1})f_2(y)
=\sum_{y\in N/N_0} f_1(y)f_2(y^{-1}x), \label{gen_prod_Hecke_alg}
\end{equation}
(where the first summation is over representatives of the right
cosets, and the second one is over representatives of the left
cosets) and the involution given by
\begin{equation}
f^*(x)=\overline{f(x^{-1})}, \label{adjoint_gen_Hecke}
\end{equation}
for $x\in N$, cf. \cite{kri}. When $x\in N$ we denote the
characteristic function of a double coset $N_0 xN_0$ variously by
$[N_0 xN_0]$ or  $[x]_{N_0}$ or simply $[x]$.

\subsection{The Hecke algebra of $(V\rtimes G, V_0\rtimes \Gamma)$}
We aim to identify conditions which will ensure that $(V\rtimes
G,V_0\rtimes\Gamma)$ is a Hecke pair and to study the Hecke algebra
of this new pair in terms of the Hecke algebras of $(G, \Gamma)$ and
$(V,V_0)$. To gain insight, we begin with a simplified but important
case.

Our convention is that for two groups $V$ and $G$ and a group
homomorphism $\alpha:G\to \Aut(V)$, the semidirect product
$V\rtimes_\alpha G$ is the cartesian product $V\times G$ endowed
with the operations
$$
(v,g)(w,h)=(v\alpha_g(w), gh)\ \ \text{and}\ \
(v,g)^{-1} =(\alpha_{g^{-1}}(v^{-1}),g^{-1})
$$
for $v,w\in V$ and $g,h\in G$. We shall omit $\alpha$ from the
notation and write  $g(v)$ for $\alpha_g(v)$ whenever convenient. We
shall also identify $V$ and $G$ with their images in $V\rtimes G$
and hence write  $vg$ instead of $(v,g)$ for a generic element of $V
\rtimes G$.

\begin{proposition}\label{medium_Hecke_pair}
Let $V$ be a group with a subgroup $V_0$, and let $\Gamma$ be a
group acting by automorphisms of $V$ in such a way that $V_0$ is
$\Gamma$-invariant. The following conditions are equivalent: \enu{i}
$(V\rtimes\Gamma,V_0\rtimes\Gamma)$ is a Hecke pair; \enu{ii}
$(V,V_0)$ is a Hecke pair and the action of $\Gamma$ on
$V_0\backslash V/V_0$ has finite orbits.

Let $\Gamma_v$ denote the stabilizer of the double coset $V_0vV_0$
for $v\in V$. If \emph{(i)} and \emph{(ii)} hold, then
$R_{V_0\rtimes\Gamma}(v)=|\Gamma/\Gamma_v|R_{V_0}(v)$, and the map
\begin{equation}
[v]_{V_0\rtimes\Gamma}\mapsto\sum_{\gamma\in\Gamma/\Gamma_v}
[\gamma(v)]_{V_0} \label{map_on_v}
\end{equation}
implements an isomorphism from
$\hh(V\rtimes\Gamma,V_0\rtimes\Gamma)$ onto $\hh(V,V_0)^\Gamma$
which respects the $R$-functions.
\end{proposition}

\begin{proof}
Note first that any right coset of $V_0\rtimes\Gamma$ is of the form
$(V_0 \rtimes \Gamma)v$ with $v\in V$. Let $v\in V$ and suppose that
\begin{equation}
(V_0\rtimes\Gamma)v(V_0\rtimes\Gamma)=\bigsqcup_{j\in
J}(V_0\rtimes\Gamma)w_j \label{decomposing_big_coset}
\end{equation}
is a disjoint decomposition of right $(V_0\rtimes\Gamma)$-cosets,
with $w_j\in V$ for all $j$ in $J$. Denote by $O_\Gamma(v)$ the orbit 
of $v$ in $V$. Intersecting the left hand side of 
\eqref{decomposing_big_coset} with $V$ gives, using the hypothesis 
$\Gamma(V_0)\subset V_0$, the set $V_0 O_\Gamma(v)V_0$. Since 
$y\in (V_0\rtimes \Gamma)z$ is equivalent to $y\in V_0z$ for all $y,z$ 
in $V$, the right hand side of \eqref{decomposing_big_coset} 
intersects $V$ in the disjoint decomposition 
$\bigsqcup_{j\in J}V_0w_j$. Thus
\begin{equation}
V_0 O_\Gamma(v)V_0=\bigsqcup_{j\in J}V_0w_j,
\label{decomp_of_orbit_DC}
\end{equation}
and the equivalence of (i) and (ii) follows.

Since $R_{V_0}(\gamma(v)) = R_{V_0}(v)$ for every $\gamma$ in
$\Gamma$, the double $V_0$-cosets contained in $V_0 O_\Gamma(v)V_0$
contain the same number of right cosets each. Thus, decomposing $V_0
O_\Gamma(v)V_0$ first into a union of double cosets, and then
decomposing each double coset into right cosets, we see that
$R_{V_0\rtimes \Gamma}(v) = |\Gamma/\Gamma_v| R_{V_0}(v)$.

The algebra $\hh(V\rtimes\Gamma,V_0\rtimes\Gamma)$ is spanned by the
characteristic functions $[v]$ of double cosets $(V_0\rtimes
\Gamma)v(V_0 \rtimes \Gamma)$ for $v\in V$.  Equation
\eqref{decomp_of_orbit_DC} implies that $[v]\vert_V$ depends only on
the $\Gamma$-orbit of $V_0vV_0$ in $V_0\backslash V/V_0$. Hence the
map in \eqref{map_on_v} is just
$\hh(V\rtimes\Gamma,V_0\rtimes\Gamma)\ni f\mapsto f|_V\in
\hh(V,V_0)^\Gamma$, and is a linear isomorphism. To see that it is
multiplicative we compute $[v]\,[w](x)$ as
$$
R_{V_0\rtimes \Gamma}((V_0\rtimes \Gamma)v\inv (V_0\rtimes
\Gamma)x\cap (V_0\rtimes\Gamma)w(V_0\rtimes \Gamma)),
$$
and we claim that this equals
\begin{equation}
\sum_{\alpha\in\Gamma/\Gamma_v,\beta\in\Gamma/\Gamma_w}
R_{V_0}(V_0\alpha(v\inv)V_0x\cap V_0\beta(w)V_0),\label{RHS}
\end{equation}
for every $v,w,x\in V$. The claim follows, upon invoking
\eqref{decomposing_big_coset} and \eqref{decomp_of_orbit_DC}, from
the decompositions
\[
(V_0\rtimes \Gamma)v\inv (V_0\rtimes \Gamma)x=
\bigcup_{\alpha\in\Gamma/\Gamma_v}(V_0\rtimes
\Gamma)V_0\alpha(v\inv)V_0x
\]
and $(V_0\rtimes \Gamma)w(V_0\rtimes \Gamma)
=\bigcup_{\beta\in\Gamma/\Gamma_w}(V_0\rtimes
\Gamma)V_0\beta(w)V_0$.
\end{proof}

As a corollary, we obtain the following generalization of
\cite[Proposition 1.7 (II.3)]{lvf}.

\begin{corollary}\label{medium_Hecke_pair_normal}
With the hypotheses of \proref{medium_Hecke_pair}, assume that $V_0$
is normal in $V$. Then there is an isomorphism
\[
\hh(V\rtimes\Gamma,V_0\rtimes\Gamma) \cong \C[V/V_0]^\Gamma,
\]
given by the restriction map $f\mapsto f\vert_V$.  Furthermore the
product in $\hh(V\rtimes\Gamma,V_0\rtimes\Gamma)$ is given by
\[
[v][w]=\sum_{\alpha\in\Gamma/\Gamma_v,\beta\in\Gamma/\Gamma_w}
|\Gamma/\Gamma_{\alpha(v)\beta(w)}|\inv [\alpha(v)\beta(w)]\ \
\text{for}\ \ v,w\in V.
\]
\end{corollary}

\begin{proof}
Since \eqref{map_on_v} is multiplicative and $V_0\trianglelefteq V$,
we obtain from  \eqref{RHS} that
\begin{equation}
[v][w](x)
=\#\{(\alpha,\beta)\in\Gamma/\Gamma_v\times\Gamma/\Gamma_w\,
|\,x=\alpha(v)\beta(w) \ \text{in}
 \ V/V_0\}.\label{computing_product}
\end{equation}
Note that $[\alpha(v)\beta(w)](x)=1$ exactly when
$\alpha(v)\beta(w)\in O_\Gamma(x)V_0$, and in this case
$$
|\Gamma/\Gamma_{\alpha(v)\beta(w)}|=|\Gamma/\Gamma_x|.
$$
We claim that
\[
\#\{(\alpha,\beta)\,|\,\alpha(v)\beta(w)=x \ \text{in}
 \ V/V_0\}=\frac{
\#\{(\alpha,\beta)\,|\,\alpha(v)\beta(w)\in O_\Gamma(x)V_0\}}
{|\Gamma/\Gamma_x|}.
\]
To prove the claim, note that $(y,z)\mapsto yz$ is a
$\Gamma$-equivariant map from $O_\Gamma(v)V_0\times O_\Gamma(w)V_0$
into~$V/V_0$, and  then the number of elements in the preimage of
$O_\Gamma(x)V_0$ under this map is the product of the number of
preimages of $xV_0\in V/V_0$ multiplied by the size of the orbit, as
needed. The result follows from the claim and \eqref
{computing_product}.
\end{proof}

We consider next the more general situation in which  there is a
group $G$ containing $\Gamma$ and acting on $V$ by an extension of
the action of $\Gamma$. To state the next result, we need to recall
from e.g. \cite[Chapter I]{kri} that two subgroups $K_1$ and $K_2$
of a group $L$ are commensurable if $K_1\cap K_2$ has finite index
in both $K_1$ and $K_2$.

\begin{proposition}
\label{HeckeSmallinsideBig} Let $G$ be a group acting by
automorphisms of a group $V$, and let $\Gamma$ be a subgroup of $G$
leaving a subgroup $V_0$ of $V$ invariant. Then $(V\rtimes
G,V_0\rtimes\Gamma)$ is a Hecke pair if and only if the following
conditions are satisfied: \enu{i} $(V,V_0)$ is a Hecke pair such
that the action of $\Gamma$ on $V_0\backslash V/V_0$ has finite
orbits and \enu{ii} $(G,\Gamma)$ is a Hecke pair such that $V_0$ and
$g(V_0)$ are commensurable subgroups of $V$ for every $g\in G$.

If \emph{(i)} and \emph{(ii)} hold, then we have an embedding
$\hh(V,V_0)^\Gamma\hookrightarrow \hh(V\rtimes G,V_0\rtimes\Gamma)$,
and
\begin{eqnarray}
R_{V_0\rtimes\Gamma}(vg) &=& R_{\Gamma}(g) \, |V_0\backslash V_0
O_{\Gamma_g}(v)g(V_0)|,\label{no_of_RC_big_Hecke}\\
L_{V_0\rtimes\Gamma}(vg) &=& L_{\Gamma}(g) \, |V_0
O_{\Gamma_g}(v)g(V_0)/g(V_0)|,\label{no_of_LC_big_Hecke}
\end{eqnarray}
where $\Gamma_g=g\Gamma g^{-1}\cap\Gamma$  and $O_{\Gamma_g}(v)$
denotes the $\Gamma_g$-orbit of $v\in V$. Furthermore,
$$
\Delta_{V_0\rtimes\Gamma}(vg)=\Delta_{V_0}(v)\Delta_\Gamma(g)\frac
{|V_0/(V_0\cap g(V_0))|}{|V_0/(V_0\cap g^{-1}(V_0))|}.
$$
\end{proposition}

\begin{proof}
We begin by proving (\ref{no_of_LC_big_Hecke}) without assuming (i)
and (ii). To simplify the notation we let $P_0: =V_0\rtimes \Gamma$.
Let $\{\gamma_i\}_i$ be representatives of left $\Gamma_g$-cosets in
$\Gamma$. Put $g_i=\gamma_i g$. Then, since $\Gamma_g$ is exactly
the set of elements $\gamma$ such that $\gamma g\Gamma=g\Gamma$,
$\{g_i\}_i$ are representatives of left $\Gamma$-cosets in $\Gamma
g\Gamma$. We claim that $P_0vgP_0$ is the union of the sets
$$
V_0\gamma_i(O_{\Gamma_g}(v))g_i(V_0)g_i\Gamma.
$$
Indeed, $P_0vgP_0$ is the union of the sets
$V_0\gamma_i\Gamma_gvgV_0\Gamma=V_0\gamma_i\Gamma_gvg(V_0)g\Gamma$.
Observe that $\gamma (g(V_0))=g(V_0)$ for $\gamma\in\Gamma_g$. Hence
$\Gamma_gvg(V_0)=O_{\Gamma_g}(v)g(V_0)\Gamma_g$, so that
$$
V_0\gamma_i\Gamma_gvg(V_0)g\Gamma
=V_0\gamma_iO_{\Gamma_g}(v)g(V_0)g\Gamma
=V_0\gamma_i(O_{\Gamma_g}(v))\gamma_i(g(V_0))\gamma_ig\Gamma,
$$
and since $\gamma_ig=g_i$, our claim is proved.

Since left $P_0$-cosets have the form $whP_0=wh(V_0)h\Gamma$, the
sets $V_0\gamma_i(O_{\Gamma_g}(v))g_i(V_0)g_i\Gamma$ do not
intersect, and
$$
L_{P_0}(vg)=\sum_i|V_0\gamma_i(O_{\Gamma_g}(v))g_i(V_0)/g_i(V_0)|.
$$
By applying the automorphisms $\gamma_i^{-1}$ we see that all
summands in the formula above are equal to
$|V_0O_{\Gamma_g}(v))g(V_0)/g(V_0)|$, and we therefore get
(\ref{no_of_LC_big_Hecke}). Then
$$
R_{P_0}(vg)=L_{P_0}(g^{-1}v^{-1})=L_{P_0}(g^{-1}(v^{-1})g^{-1})
=L_\Gamma(g^{-1})\,|V_0O_{\Gamma_{g^{-1}}}(g^{-1}(v^{-1}))g^{-1}(V_0)
/g^{-1}(V_0)|.
$$
Applying the automorphism $g$ we get that the last expression equals
$$
L_\Gamma(g^{-1})\,|g(V_0)O_{g\Gamma_{g^{-1}}g^{-1}}(v^{-1})V_0
/V_0|=R_\Gamma(g)|g(V_0)O_{\Gamma_g}(v^{-1})V_0 /V_0|.
$$
Finally, by applying the anti-automorphism $V\ni w\mapsto w^{-1}$,
we get
$$
R_{P_0}(vg)=R_\Gamma(g)|V_0\backslash V_0O_{\Gamma_g}(v)g(V_0)|,
$$
and (\ref{no_of_RC_big_Hecke}) is also proved.

It is now immediate that $(V\rtimes G,V_0\rtimes\Gamma)$ is a Hecke
pair if and only if (i) and (ii) hold. That we have an embedding
$\hh(V,V_0)^\Gamma\hookrightarrow \hh(V\rtimes G,V_0\rtimes\Gamma)$
follows from \proref{medium_Hecke_pair}.

It remains to prove the formula for the modular function (note that
it will also follow from \secref{subsection_Schlichting_completion}
below). We have
$$
\Delta_{P_0}(vg)=\Delta_\Gamma(g)\frac{|V_0
O_{\Gamma_g}(v)g(V_0)/g(V_0)|}{|V_0\backslash
V_0 O_{\Gamma_g}(v)g(V_0)|}.
$$
Since the numbers $|V_0 wg(V_0)/g(V_0)|$ and $|V_0\backslash
V_0wg(V_0)|$ are independent of $w\in O_{\Gamma_g}(v)$, we conclude
that $\Delta_{P_0}(vg)$ is equal to
$$
\Delta_\Gamma(g)\frac{|V_0vg(V_0)/g(V_0)|} {|V_0\backslash
V_0vg(V_0)|}= \Delta_\Gamma(g)\frac{|v^{-1}V_0vg(V_0)/g(V_0)|}
{|v^{-1}V_0v\backslash v^{-1}V_0vg(V_0)|}
=\Delta_\Gamma(g)\frac{|v^{-1}V_0v/(v^{-1}V_0v\cap g(V_0))|}
{|g(V_0)/(v^{-1}V_0v\cap g(V_0))|}.
$$
If $K_1$, $K_2$ and $K_3$ are commensurable groups then
$$
\frac{|K_1/(K_1\cap K_2)|}{|K_2/(K_1\cap K_2)|} =\frac{|K_1/(K_1\cap
K_3)|}{|K_3/(K_1\cap K_3)|}\, \frac{|K_3/(K_2\cap
K_3)|}{|K_2/(K_2\cap K_3)|}.
$$
Applying this to $K_1=v^{-1}V_0v$, $K_2=g(V_0)$ and $K_3=V_0$ we get
\begin{align*}
\Delta_{P_0}(vg) &=\Delta_\Gamma(g)\frac{|v^{-1}V_0v/(v^{-1}V_0v\cap
V_0)|}{|V_0/(v^{-1}V_0v\cap V_0)|}\, \frac{|V_0/(g(V_0)\cap
V_0)|}{|g(V_0)/(g(V_0)\cap
V_0)|}\\
&=\Delta_\Gamma(g)\Delta_{V_0}(v)\frac{|V_0/(g(V_0)\cap
V_0)|}{|g(V_0)/(g(V_0)\cap V_0)|},
\end{align*}
which completes the proof of the proposition.
\end{proof}

\subsection{A Hecke pair from the inclusion $\mz \subset \mq$}
\label{favouritepair}

We now explain our motivating example, which is the  $2\times 2$
matrix analogue of the Hecke pair $(P_{\Q}^+, P_{\Z}^+)$ of Bost and
Connes \cite{bos-con}, and arises from the following setup: let
$V:=\mq$ and $V_0:=\mz$ be the additive groups in  the ring
inclusion $\mq\supset \mz$, and  let $G $ and $\Gamma $ be their
orientation preserving groups of invertible elements acting by
multiplication on the right. Specifically,
\begin{eqnarray*}
G &=& \glq:=\{g\in \mq\mid \det(g)>0\},\\
\Gamma &=& \slz :=\{m\in \mz\mid \det(m)=1\},
\end{eqnarray*}
and $\alpha_g(m)=mg^{-1}$ for $g\in \glq$ and $m\in \mq$. Since
$(0,g)(m,1)(0,g^{-1})=(\alpha_g(m),1)$ is compatible with the matrix
multiplication
$$
\matr{0}{g}\matr{m}{1}\matr{0}{g^{-1}}=\matr{mg^{-1}}{1},
$$
we can view the element $(\alpha_g(m), g)$ of $V\rtimes G$ as the
$4\times 4$ matrix $\matr{m}{g}$.

To argue that $(V\rtimes G,V_0\rtimes\Gamma)$ is a Hecke pair does
not require any computations, because it can be deduced from the
following general result:  if~$\hh$ is a linear algebraic group
defined over~$\Q$, and $H=\hh(\Q)$ and $H_0=\hh(\Z)$ are the groups
of rational and integral points of $\hh$, then $(H,H_0)$ is a Hecke
pair, see e.g. \cite[Corollary~1 of~Proposition~4.1]{pr}. However,
to compute the $R$-function we  need the elementary considerations
of Proposition~\ref{HeckeSmallinsideBig} anyway, so at this point we
prefer to rely on that proposition.

It is well known that  $(G, \Gamma)$ is a Hecke pair, in fact
$\hecke{G}{\Gamma}$ is the classical Hecke algebra generated by the
Hecke operators, see e.g. \cite[Chapter IV]{kri} for details. The
remaining assumptions of Proposition~\ref{HeckeSmallinsideBig} are
easily verified.   E.g. to show that $\Gamma$-orbits in $V/V_0$ are
finite, take $m\in \mq$ and choose $d\in \N^*$ such that $dm \in
\mz$. Then $dm \slz\subset\mz$, and thus the $\slz$-orbit of $m +\mz$
in $\mq/\mz$ is finite, and in fact has at most $ |\mz d^{-1}/\mz| =
d^4$ elements. Thus $(V\rtimes G,V_0\rtimes\Gamma)$ is a Hecke pair,
and since $(G,\Gamma)$ is unimodular, we have
$$
\Delta_{V_0\rtimes\Gamma}\matr{m}{g}=\frac {|\mz/(\mz\cap \mz
g^{-1})|}{|\mz/(\mz\cap \mz g)|}=\det(g)^{-2}.
$$
We have therefore obtained the following.

\begin{proposition}\label{MainHecke}
The inclusion $ \mz \rtimes \slz \subset  \mq \rtimes \glq$ is a
Hecke pair. Furthermore the associated modular function is given by
$\Delta\matr{m}{g}=\det(g)^{-2}.$
\end{proposition}

\subsection{Hecke algebras from subsemigroups}
\label{subsec_semi_Hecke} Suppose that $(N, N_0)$ is a Hecke pair
and $S$ is a subsemigroup of $N$  containing $N_0$. The set $\hh(S,
N_0)$ of functions in $\hh(N, N_0)$ that are supported on the double
cosets of elements of $S$ forms an  algebra, with convolution
formula \eqref{gen_prod_Hecke_alg} restricted to a summation over
$y$ in $S/N_0$ but without adjoints. It is easy to see that $\hh(S,
N_0)$ is a subalgebra of $\hh(N, N_0)$, cf. \cite[Lemma I.4.9]{kri}.

\begin{example}
The main example of this situation is again classical: for the Hecke
pair $(G, \Gamma)=(\glq, \slz)$ of \secref{favouritepair}, we take
$S$ to be the subsemigroup of integer matrices with positive
determinant, $\mzp :=\{g\in \mz\mid \det(g)>0\}$; the Hecke algebra
$\hh(S, \Gamma)$ is described in \cite[Chapter IV]{kri}.
\end{example}

We have seen in \proref{HeckeSmallinsideBig} that $\hh (V \rtimes G,
V_0\rtimes \Gamma)$ contains a copy of the fixed point algebra of
$\hh(V,V_0)$ under the action of $\Gamma$, so it is natural to
wonder whether $\hh(G,\Gamma)$ plays any role in $\hh (V \rtimes G,
V_0\rtimes \Gamma)$. The following proposition shows that there is
an actual algebraic embedding, but only of the nonselfadjoint Hecke
subalgebra of $\hh(G,\Gamma)$ corresponding to a convenient
subsemigroup of $G$.

\begin{proposition} \label{smallHecke}
Suppose that $(V\rtimes G, V_0\rtimes \Gamma)$ is a Hecke pair as in
Proposition~\ref{HeckeSmallinsideBig} and let $S= \{ s\in G\mid
V_0\subset s(V_0)\}$. Then $S$ is a semigroup containing $\Gamma$
and  the map
\[
\iota \colon  [s]_\Gamma \in \hh(S, \Gamma) \mapsto [s]_{V_0\rtimes
\Gamma} \in \hh(V\rtimes G,V_0\rtimes \Gamma)
\] extends by linearity to an injective homomorphism.
\end{proposition}

\bp Clearly $S$ is a semigroup. To prove the second assertion it is
convenient to use the notation $P =V\rtimes G $ and $P_0 = V_0
\rtimes \Gamma$. The set $SV_0$ is a subsemigroup of $V\rtimes G$
containing $V_0\rtimes\Gamma$, so $\hecke{SV_0}{P_0}$ is a
subalgebra of $\hecke{P}{P_0}$. Furthermore, the identity map on $S$
gives a bijection between $\Gamma\backslash S/\Gamma$ and
$P_0\backslash SV_0/P_0$. It follows that $\iota$ extends to a
bijective linear map from $\hecke{S}{\Gamma}$ onto
$\hecke{SV_0}{P_0}$. To see that this is an isomorphism of algebras,
we have to show that
$$
R_{P_0}(P_0t^{-1}P_0r\cap P_0sP_0)=R_\Gamma(\Gamma t^{-1}\Gamma
r\cap \Gamma s\Gamma)
$$
for $r,s,t\in S$. Since $P_0t^{-1}P_0r= V_0 \Gamma t^{-1} \Gamma r$,
we have $P_0t^{-1}P_0r\cap P_0sP_0= V_0(\Gamma t^{-1} \Gamma r\cap
\Gamma s \Gamma)$, as required.  \ep

It is easy to see that the embedding $\iota$ of \proref{smallHecke}
restricts to the Hecke algebra of any subsemigroup of $S$ containing
$\Gamma$. But the embedding does not extend in general to a
$*$-preserving homomorphism from  $\hh(G,\Gamma)$ to $\hh(V\rtimes
G,V_0\rtimes \Gamma)$; in particular, this happens in our main
example, essentially because the classical Hecke algebra
$\hh(\glq,\slz)$  is commutative, but double cosets of $\mz \rtimes
\slz$ do not $*$-commute in general.

\begin{remark}\label{another_homo}
Under the assumptions of \proref{smallHecke}, if $\phi\colon
S\to\R^*$ is a homomorphism whose kernel contains $\Gamma$, then by
essentially the same proof, $[s]_\Gamma\mapsto\phi(s)[s]_{V_0\rtimes
\Gamma}$ is also a homomorphism from $\hh(S, \Gamma)$ into
$\hh(V\rtimes G, V_0\rtimes \Gamma)$. We shall see in
\secref{subsection_corner_rep} that the choice $\phi(s):=
[s(V_0):V_0]^{-1/2}$ instead of $\phi(s)$=1 gives rise to a more
natural homomorphism of Hecke algebras.
\end{remark}

\bigskip

\section{Completions of Hecke pairs}\label{completions}

If $N$ is a locally compact group and $N_0$ is a compact open
subgroup then $(N,N_0)$ forms a Hecke pair. It is known that any
Hecke pair gives rise to a pair of this form with isomorphic Hecke
algebra, and we shall review the known facts, but also describe
slightly more general topological Hecke pairs that turn out to be
useful in the context of semidirect product groups.

\subsection{Completions of Hecke pairs and their C$^*$-algebras}
\label{stuff_compl} Let $(N,N_0)$ be a Hecke pair.

\begin{definition} \label{Not_Schlichting}
A pair $(\bar N, \bar N_0)$ and a group homomorphism $\rho\colon
N\to \bar N$ form a \emph{completion} of a Hecke pair $(N, N_0)$ if
$\bar N$ is a locally compact group, $\bar N_0$ is a compact open
subgroup of $\bar N$, $\rho(N)$ is dense in $\bar N$, $\rho(N_0)$ is
dense in $\bar N_0$, and $\rho^{-1}(\bar N_0)=N_0$.
\end{definition}

Completions always exist, and in fact one can produce a canonical
one, called the \emph{Schlichting completion}~\cite{tza,klq}, see
also  \cite{glo-wil,larr} for slightly different approaches. We
recall from~\cite{tza,klq} that a Hecke pair $(N, N_0)$  is called
\emph{reduced} if $ \bigcap_{x\in N}xN_0x^{-1}=\{e\}$ or,
equivalently, if $N_0$ contains no nontrivial normal subgroup of
$N$. The \emph{Hecke topology} on $N$ is the topology determined by
a neighbourhood subbase at $e$ consisting of the sets $xN_0x^{-1}$
for $x\in N$. Then the Schlichting completion of $(N, N_0)$ is a
Hecke pair $(\bar N, \bar N_0)$ together with a homomorphism
$\phi\colon N\to \bar N$ characterized uniquely  by the properties
that $\bar{N}$ is  a locally compact totally disconnected group,
$\bar N_0$ is a compact open subgroup, $(\bar{N}, \bar N_0)$ is
reduced, $\phi(N)$ is dense in $\bar{N}$, and $\phi^{-1}(\bar
N_0)=N_0$. Furthermore, if $(N,N_0)$ is reduced then $\phi$ is a
homeomorphism of $N$ with its Hecke topology onto its image inside
$\bar N$. The Schlichting completion is universal in the following
sense.

\begin{proposition}\label{Hecke_alg_same_for_non_Schlichting}
Let $(N, N_0)$ be a Hecke pair, and suppose that $(\tilde{N},
\tilde{N_0})$ and $\rho\colon N\to \tilde{N}$ form a completion of
$(N, N_0)$. Assume also that $(\bar{N}, \bar N_0)$ together with a
homomorphism $\phi\colon N\to \bar{N}$ is the Schlichting completion
of $(N,N_0)$. Then there exists a unique continuous homomorphism
$\tilde\phi\colon\tilde N\to\bar{N}$ such that
$\tilde\phi\circ\rho=\phi$. The homomorphism $\tilde\phi$ is onto,
and it is an isomorphism if and only if the pair $(\tilde N,\tilde
N_0)$ is reduced.
\end{proposition}

\bp For the proof note that we can identify $N/N_0$ with $\tilde
N/\tilde N_0$, and then $\bar N$ is the Schlichting completion of
$\tilde{N}$  by construction, see e.g.~\cite{klq}. \ep

Let now $(\bar N,\bar N_0)$ together with $\rho\colon N\to\bar N$ be
a completion of $(N,N_0)$. Denote by $\mu$ the left Haar measure on
$\bar N$ such that $\mu(\bar N_0)=1$, and by $\Delta=\Delta_{\bar
N}$ the modular function of $\bar N$. Note that by \cite[Lemma
1]{schl} the composition of $\Delta_{\bar N}$ with $\rho$ coincides
with $\Delta_{N_0}$ from \secref{subsection_hecke_alg}.

The space $C_c(\bar N)$ of compactly supported continuous functions
on $\bar N$ is a $*$-algebra with the convolution product
\begin{equation}
(f_1\ast f_2)(x)=\int_{\bar N} f_1(y)f_2(y^{-1}x)d\mu(y) =\int_{\bar N}
\Delta(y^{-1})f_1(xy^{-1})f_2(y)d\mu(y)\label{def_conv_on_Cc}
\end{equation}
and involution
\begin{equation}
f^*(x)=\Delta(x^{-1})\overline{f(x^{-1})}.\label{def_inv_on_Cc}
\end{equation}

Denote by $p$ the characteristic function~$\31_{\bar N_0}$ of $\bar
N_0$, which is a  self-adjoint projection in $C_c(\bar{N})$. Then we
have the following result similar to the situation of Schlichting
completions in~\cite{tza, klq} (note that in \cite{klq} the
definition of the involution has an extra factor
$\Delta_{N_0}(x^{-1})$).

\begin{lemma}\label{Hecke_alg_as_*_corner}
The homomorphism $\Psi\colon \hh(N,N_0)\to C_c(\bar{N})$,
$\Psi([N_0xN_0]):= \Delta_{N_0}^{-1/2}(x)[\bar N_0\rho(x)\bar N_0]$, 
is a $*$-isomorphism from $\hh(N,N_0)$ onto the $*$-algebra
$pC_c(\bar{N})p$ with operations
\eqref{def_conv_on_Cc}--\eqref{def_inv_on_Cc}.
\end{lemma}

Suppose further that $\pi\colon \bar N\to U(H)$ is a unitary
representation of $\bar N$ on a Hilbert space $H$. The integrated
form $\pi_*$ of $\pi$ carries $p$ into the projection in $B(H)$ onto
the space $\{\xi\in H\mid \pi(n)\xi=\xi, \forall\ n \in  \bar N_0\}$
of $\bar N_0$-fixed vectors. Hence by \lemref{Hecke_alg_as_*_corner}
we obtain  a $*$-representation
$$
\pi_*\circ \Psi\colon \hh(N,N_0)\to B(\pi_*(p)H).
$$
In particular, when $\pi$ is the regular representation of
$\bar{N}$, we get a representation $\pi_*\circ \Psi$ of $\hh(N,
N_0)$ on the space of $\bar N_0$-fixed vectors in $L^2(\bar{N})$.
The map $U\colon\ell^2(N_0\backslash N)\to L^2(\bar{N})$ defined by
$$
(U\xi)(n):=\Delta_{\bar N}(n)^{-1/2}\xi(x),
$$
where $x$ is such that $n\in \bar N_0\rho(x)$, is an isometry of
$\ell^2(N_0\backslash N)$ onto the space of $\bar N_0$-fixed vectors
in~$L^2(\bar{N})$. It follows that $U^*(\pi_*\circ \Psi)(\cdot)U$ is
a $*$-representation of $\hh(N, N_0)$ on $\ell^2(N_0\backslash N)$.
We claim that this is the familiar \emph{regular representation}
$\lambda$ of $\hh(N, N_0)$ from \cite[Proposition 3]{bos-con}, which
satisfies
\begin{equation}
(\lambda(f)\xi)(x)=\sum_{y\in N_0 \backslash N} f(xy^{-1})\xi(y)
=\sum_{y\in N/N_0} f(y)\xi(y^{-1}x).
\label{def_originallrr_Hecke_alg}
\end{equation}
Indeed, denoting by $s\colon\bar N\to N$ a map such that $n\in\bar
N_0\rho(s(n))$ for $n\in\bar N$, we compute
\begin{equation}\label{how_to_get_lrr}
\begin{split}
((\pi_*\circ\Psi)(f)U\xi)(m)
&=\int_{\bar{N}}\Psi(f)(n)(U\xi)(n^{-1}m)d\mu(n)\\
&=\Delta(m)^{-1/2}\sum_{y\in N/N_0}\int_{\rho(y)\bar N_0}
f(s(n))\xi(s(n^{-1}m))d\mu(n)\\
&=\Delta(m)^{-1/2}\sum_{y\in N/N_0}\int_{\bar N_0}
f(y)\xi(y^{-1}s(m))d\mu(n)\\
&=\Delta(m)^{-1/2}\sum_{y\in N/N_0} f(y)\xi(y^{-1}s(m)).
\end{split}
\end{equation}

The \emph{reduced Hecke C$^*$-algebra} $\redheck{N}{N_0}$ is the
C$^*$-algebra generated by the image of $\hh(N, N_0)$ in
$B(\ell^2(N_0\backslash N))$ under the representation $\lambda$.
 The isomorphism $\Psi\colon\hecke{N}{N_0}\to
pC_c(\bar{N})p$ extends to an isomorphism of $\redheck{N}{N_0}$ onto
$p\red{\bar{N}}p$, which we shall still denote by $\Psi$.

Notice that if $S$ is a subsemigroup of $N$ containing $N_0$ then
the formula \eqref{def_originallrr_Hecke_alg} makes sense when $f\in
\hh(S, N_0)$, $\xi\in \ell^2(N_0\backslash S)$ and the summation is
over $y\in N_0\backslash S$. Hence we can define the regular
representation of $\hh(S, N_0)$ on $\ell^2(N_0\backslash S)$ by
\begin{equation}
(\lambda(f)\xi)(s)=\sum_{t\in N_0\backslash S}f(st^{-1})\xi(t)
\text{ for }s\in S. \label{def_regularr_semiHecke}
\end{equation}

\subsection{The Schlichting completion of $({\mq}\rtimes {\glq},
{\mz}\rtimes{\slz})$} \label{scompletion} We first recall
from~\cite[Section 5.1]{pr} a standard way of getting completions of
Hecke pairs arising from algebraic groups (see also~\cite[Section
4]{tza}).

Let $\rr$ denote the ring $
\prod_ {p\in \mathcal{P}}\Z_p$,  where $\mathcal{P}$ is the set of
prime numbers. The additive group of $\rr$ is a compact group. The
ring of finite adeles is, by definition,
$$
\af=\prod_{p\in \mathcal{P}} (\Q_p:\Z_p),
$$
and as an additive group has the locally compact totally
disconnected topology of a restricted product in which  $\rr$ is a
compact open subgroup. Then $\Q$ embeds diagonally into $\af$, and
$\Q\cap\rr=\Z$. Consider now the group $\gln(\af)$ where  the
topology is defined by the embedding $\gln(\af)\ni g\mapsto
(g,\det(g)^{-1})\in\mn(\af)\times\af$. Equivalently, $\gln(\af)$ is
the restricted topological product of the groups $\gln(\Q_p)$ with
respect to the subgroups $\gln(\Z_p)$, so that $\gln(\rr)$ is a
compact open subgroup of $\gln(\af)$.

Now if $\hh\subset\gln$ is a linear algebraic group defined over
$\Q$, and $H=\hh(\Q)$, $H_0=\hh(\Z)$, we have an embedding
$H\hookrightarrow\gln(\af)$, which together with the groups $\tilde
H:=\overline{\hh(\Q)}$ and $\tilde H_0:=\tilde
H\cap\gln(\rr)=\overline{\hh(\Z)}$ forms a completion of $(H,H_0)$.

Returning to our Hecke pair from \secref{favouritepair}, we have an
algebraic subgroup $\hh:=\mtwo\rtimes\gl$ in $\glf$, and then
$\mq\rtimes\glq$ is a subgroup of index $2$ in $\hh(\Q)$. Since its
closure in $\overline{\hh(\Q)}\subset \glf(\af)$ will then be open,
to find a completion it suffices to find the closure of
$$
\matr{\mq}{\glq}
$$
in $\glf(\af)$. Since the group $\mq$ is dense in $\ma$, it remains
to find the closure of $\glq$ in $\gl(\af)$. This requires more care
than the one dimensional case, in which   $\Q^*_+$ is discrete in
$\glo(\af)$ because $\Q^*_+\cap\rr^*=\{1\}$; now by the strong
approximation theorem \cite[Theorem 7.12]{pr}, the closure of $\slz$
inside $\gla$ is $\slr$, the compact subgroup of matrices in
$\G(\rr)$ which have determinant $1$, see also
\cite[Lemma~1.38]{shi} for an elementary proof of our particular
case.

\begin{lemma}\label{glqdenseinglap}
The closure of $\glq$ in $\gla$ is the set
$$
\glap:=\{g\in \gla\mid \det(g)\in \Q_+^*\}.
$$
Moreover,
\begin{equation}
\glap = \glq \slr. \label{decomp_of_glap}
\end{equation}
\end{lemma}

Note that the factorization in \eqref{decomp_of_glap} is not unique.
The factors are only determined up to an element of $\slz =
\glq\cap\slr$.

\begin{proof}
Since the determinant function is continuous on $\gla$ and $\Q^*$ is
closed in $\glo(\af)$, the subgroup $\glap$ is closed, and hence is
a totally disconnected, locally compact group in its own right. In
particular, the closure of $\glq$ is contained in $\glap$.

Clearly $\glq  \slr \subset \glap$. In order to prove the reverse
inclusion, let $g\in \glap$. We need to show that $g = g_0r$ for
matrices $g_0 \in \glq$ and $r\in \slr$.  Since $g\rr^2 $ is a
lattice in~$\af^2$, that is, a compact open $\rr$-submodule, by
\cite[Theorem V.2]{weil} there exists $g_0 \in \glq$ such that
$g\rr^2 \cap \Q^2 = g_0\Z^2$. Thus $g\rr^2 = g_0 \rr^2$, which
implies that $g_0\inv g \in \glr$. Since $g$ has positive
determinant, $\det (g_0\inv g) \in \Q_+^* \cap \rr^* = \{1\}$, and
hence $g = g_0 (g_0\inv g) \in \glq\slr$, as required. As $\slz$ is
dense in $\slr$, this in particular implies that $\glq$ is dense in
$\glap$.
\end{proof}

We therefore obtain the following result.

\begin{proposition}\label{completion_our_pair}
The pair $(\ma\rtimes \glap, \mr\rtimes \slr)$ is the Schlichting
completion of $({\mq}\rtimes {\glq}, {\mz}\rtimes {\slz})$.
\end{proposition}

\bp That $(\ma\rtimes \glap, \mr\rtimes \slr)$, with the obvious
embedding of ${\mq}\rtimes {\glq}$ into $\ma\rtimes \glap$, form a
completion follows from the preceding lemma and the discussion about
algebraic groups. To see that it  is in fact the Schlichting
completion, notice that by
Proposition~\ref{Hecke_alg_same_for_non_Schlichting} we just have to
check that $(\ma\rtimes \glap, \mr\rtimes \slr)$ is a reduced pair,
which can be done using the following simple sufficient condition.
\ep

\begin{lemma} \label{suff_reduced}
A Hecke pair $(V\rtimes G,V_0\rtimes\Gamma)$ is reduced if the
following two conditions are satisfied: \enu{i} $\bigcap_{g\in
G}g(V_0)=\{e\}$; \enu{ii} if $\gamma(v)=v$ for some
$\gamma\in\bigcap_{g\in G}g\Gamma g^{-1}$ and all $v\in V$ then
$\gamma=e$.
\end{lemma}

\bp Assume $v\gamma$ lies in a normal subgroup of $V\rtimes G$ which
is contained in $V_0\rtimes\Gamma$. Then $v\gamma$ belongs to
$gV_0\Gamma g^{-1}=g(V_0)g\Gamma g^{-1}$ for all $g\in G$, which
shows that $v\in\bigcap_{g\in G}g(V_0)$ and $\gamma\in\bigcap_{g\in
G}g\Gamma g^{-1}$. Hence $v=e$ by condition (i). But then $\gamma
wg=\gamma(w)\gamma g$ belongs to $wgV_0\Gamma=wg(V_0)g\Gamma$ for
all $w\in V$, so that in particular $\gamma(w)\in\bigcap_{g\in
G}wg(V_0)=\{w\}$. Hence $\gamma=e$ by condition (ii). \ep

\subsection{Completions of $(V\rtimes G, V_0\rtimes \Gamma)$}
\label{subsection_Schlichting_completion} If $(V\rtimes G,
V_0\rtimes\Gamma)$ is a Hecke pair, there is no apparent reason for
its Schlichting completion to be again a semidirect product,
although this is the case in our main example. However, it turns out
that if we  do not require a reduced pair as a completion, we can
get one that is a semidirect product.

\begin{remark}\label{Braconnier_top}
Let $\bar V$ be a locally compact group. Let $\HA(\bar V)$ be the
group of homeomorphic group automorphisms of $\bar V$. It is a
topological group for the topology that has as a neighbourhood base
at the identity the sets
$$
\{\beta\in \HA(\bar V)\mid \beta(u)\in Uu \text{ and
}\beta^{-1}(u)\in Uu\text{ for all } u\in K\},
$$
for all compact subsets $K$ of $\bar V$ and all neighbourhoods $U$
of $e$ in $\bar V$, see \cite[(26.5)]{hew-ros}

If $\bar G$ is another locally compact group with an action $\alpha$
on $\bar V$, then this action is continuous, in the sense that the
map $\bar G\times \bar V \to \bar V$ is continuous, if and only if
the homomorphism $\alpha \colon\bar G \to \HA(\bar V)$ is
continuous. In this case $\bar V\rtimes_\alpha \bar G$ is a locally
compact group with the topology inherited from the product topology
on $\bar V\times \bar G.$
\end{remark}

\begin{proposition} \label{sdp_completion_universal}
Suppose that $(\overline{V\rtimes G}, \overline{V_0\rtimes \Gamma})$
and $\rho\colon V\rtimes G\to \overline{V\rtimes G}$ form a
completion of $(V\rtimes G,V_0\rtimes\Gamma)$, and  for each subset
$X\subset V\rtimes G$ let $X^\sharp$ denote the closure of $\rho(X)$
in $\overline{V\rtimes G}$.  Then there is a continuous action of
${G}^\sharp$ by homeomorphic automorphisms of ${V}^\sharp$ such that
$({V}^\sharp\rtimes {G}^\sharp, {V_0}^\sharp\rtimes
{\Gamma}^\sharp)$ and the homomorphism $\iota\colon (v,g)\mapsto
(\rho(v),\rho(g))$ from $V\rtimes G$ into ${V}^\sharp\rtimes
{G}^\sharp$ form a completion of $(V\rtimes G, V_0\rtimes \Gamma)$.
\end{proposition}

\begin{proof}
Since $\rho(G)$ normalizes $\rho(V)$, it normalizes $V^\sharp$.
Hence $G^\sharp$ normalizes $V^\sharp$. We therefore have a
continuous action of $G^\sharp$ on $V^\sharp$ by conjugation.

Since both $V_0^\sharp$ and $\Gamma^\sharp$ are compact,
$V_0^\sharp\Gamma^\sharp$ is a closed subgroup of
$\overline{V_0\rtimes \Gamma}$, and hence
$V_0^\sharp\Gamma^\sharp=\overline{V_0\rtimes \Gamma}$. We claim
that $V^\sharp \cap (V_0^\sharp\Gamma^\sharp)=V_0^\sharp$. Indeed,
take $v$ in $V^\sharp \cap (V_0^\sharp\Gamma^\sharp)$, and write $v=
\lim_i \rho(v_i)$ for $v_i\in V$. The set $V_0^\sharp\Gamma^\sharp$
being open, it eventually contains $\rho(v_i)$, and so $ v_i\in
V\cap \rho^{-1}(\overline{V_0\rtimes \Gamma})=V\cap (V_0\rtimes
\Gamma)=V_0$, proving that $v\in V_0^\sharp$. The claim implies that
$V_0^\sharp$ is open in $V^\sharp$, and that
$$
V\cap\rho^{-1}(V_0^\sharp)\subset V\cap
\rho^{-1}(\overline{V_0\rtimes \Gamma})=V\cap (V_0\rtimes \Gamma)=
V_0.
$$
Hence $V\cap\rho^{-1}(V_0^\sharp)=V_0$, and thus the pair
$(V^\sharp, V_0^\sharp)$ and $\rho\vert_{V}$ form a completion of
$(V,V_0)$. Similarly, $(G^\sharp, \Gamma^\sharp)$ and
$\rho\vert_{G}$ form a completion of $(G, \Gamma)$, and then the map
$\iota=(\rho\vert_V, \rho\vert_G)$ satisfies the claim of the
proposition.
\end{proof}

The next theorem shows that  it is possible to produce completions
of semidirect product pairs from completions of  the component
pairs. To simplify the notation we shall only consider reduced pairs
and identify a group with its image in the  Schlichting completion,
suppressing the corresponding injective homomorphism.

\begin{theorem}\label{SemidirectCompletions}
Let $(V\rtimes G, V_0 \rtimes \Gamma)$ be a Hecke pair as in
\proref{HeckeSmallinsideBig}, and suppose that it is reduced. Denote
by $\bar{V}$ and $\bar V_0$ the closures of $V$ and $V_0$ in the
Schlichting completion of $(V\rtimes G, V_0)$, and suppose that
$(\tilde{G}, \tilde{\Gamma})$ and $\tilde{\rho}\colon G\to
\tilde{G}$ form a completion of $(G, \Gamma)$. \enu{i} The map
$\rho(g):=(\Ad g, \tilde{\rho}(g))$ is a homomorphism from $G$ into
$\HA(\bar{V})\rtimes \tilde{G}$, and the closures
$\bar{G}:=\overline{\rho(G)}$ and
$\bar{\Gamma}:=\overline{\rho(\Gamma)}$ satisfy that
$(\bar{V}\rtimes \bar{G}, \bar V_0\rtimes \bar{\Gamma})$ is a
completion of $(V\rtimes G,V_0\rtimes \Gamma) $ together with the
map $\iota\colon V\rtimes G \to\bar{V}\rtimes \bar{G}$ given by
$\iota\colon (v,g)\mapsto (v, \rho(g))$. \enu{ii} Let
$(\overline{V\rtimes G}, \overline{V_0\rtimes \Gamma})$ be the
Schlichting completion of $(V\rtimes G, V_0\rtimes \Gamma)$. If
$(\tilde{G}, \tilde{\Gamma})$ is the Schlichting completion of $(G,
\Gamma)$ and for every $v$ in $V$ there is a finite set
$\{g_1,\dots,g_n\}$ in $G$ such that
\begin{equation}
\bigcap_{i=1}^ng_i(V_0)\subset vV_0v^{-1},
\label{cond_for_Schlichting}
\end{equation}
then $\iota$ extends to a topological isomorphism of
$\overline{V\rtimes G}$ onto $\bar{V}\rtimes {\bar{G}}$ and of
$\overline{V_0\rtimes \Gamma}$ onto $\bar V_0\rtimes \bar{\Gamma}$.
\label{getting_Schlichting_semidirect}
\end{theorem}

\begin{proof}
We claim that the closure of the image of $\Gamma$ in $\HA(\bar{V})$
under the map $\Ad$ is compact. To see this  note first that
$(V\rtimes G, V_0\rtimes \Gamma)$ being reduced implies on one hand
that $(V\rtimes G, V_0)$ is reduced, and on the other that
$\overline{V\rtimes G}$  has the Hecke topology, given by a
neighbourhood subbase at $e$ consisting of sets of the form
$(v,g)(V_0\rtimes \Gamma)(v,g)^{-1}$ for $(v,g)\in V\rtimes G$. The
relative topology on $V$ has the sets of the form $vg(V_0)v^{-1}$ as
elements of the subbase, and these are precisely the sets defining
the Hecke topology on $V\rtimes G$ for the pair $(V\rtimes G, V_0)$.
Hence~$\bar{V}$ can be considered as a closed normal subgroup of
$\overline{V\rtimes G}$.

The closure of $V\rtimes G$ in the Hecke topology coming from
$(V\rtimes G, V_0)$ acts on $\bar{V}$ by conjugation, and this
action drops to an action of $G$ on $\bar{V}$. Since the closure
of~$\Gamma$ in~$\overline{V\rtimes G}$ is compact, the closure of
$\Gamma$ in $\HA(\bar{V})$ is compact, so our claim is proved.

We claim next that $(\bar{G}, \bar{\Gamma})$ and the map $\rho\colon
g\mapsto (\Ad g,\tilde{\rho}(g))$ form a completion of $(G,
\Gamma)$. That $\bar{\Gamma}$ is compact in the product
$\HA(\bar{V})\times \tilde{G}$ follows because this is true in each
component. We have $ {\rho}^{-1}(\bar{\Gamma})\subset
\tilde{\rho}^{-1}(\tilde{\Gamma})=\Gamma, $ hence
${\rho}^{-1}(\bar{\Gamma})=\Gamma$, as required. To see that
$\bar{\Gamma}$ is open, let $\pi_2\colon \bar{G}\to \tilde{G}$ be
the projection map, and take $x$ in $\pi_2^{-1}(\tilde{\Gamma})$. We
can approximate $x$ by an element of the form $\rho(g)$ for $g\in
G$, and then $\tilde{\rho}(g)=\pi_2(\rho(g)) \in \tilde{\Gamma}$,
showing that $g\in \tilde{\rho}^{-1}(\tilde{\Gamma})= \Gamma$. Hence
$x\in \bar{\Gamma}$, and so $\bar{\Gamma}$ is
$\pi_2^{-1}(\tilde{\Gamma})$, and is therefore open, proving the
claim.

Since $\bar{G}$ has a continuous action on $\bar{V}$ by
construction, and since this action gives by restriction to
$\bar{\Gamma}$ an action on $\bar V_0$, the pair $(\bar{V}\rtimes
\bar{G}, \bar V_0\rtimes \bar{\Gamma})$ and the map $\iota$ satisfy
the conditions of Definition~\ref{Not_Schlichting}, and we have
proved (i).

To prove (ii) it suffices by
Proposition~\ref{Hecke_alg_same_for_non_Schlichting} to ensure
that $(\bar{V}\rtimes \bar{G}, \bar V_0\rtimes \bar{\Gamma})$ is
reduced. We shall check that the conditions of
Lemma~\ref{suff_reduced} are satisfied.

Since the Schlichting completion of $(V\rtimes G,V_0)$ is by
definition reduced, we have
$$
\bigcap_{g\in G,v\in V}g(v\bar V_0v^{-1})=\bigcap_{g\in G,v\in
V}gv\bar V_0v^{-1}g^{-1}=\{e\}.
$$
By (\ref{cond_for_Schlichting}) the left hand side coincides with
$\bigcap_{g\in G}g(\bar V_0)$, and thus assumption (i) of
Lemma~\ref{suff_reduced} is satisfied.

Assume now that $\gamma\in\bigcap_{g\in
G}\rho(g)\bar{\Gamma}\rho(g)^{-1}$ acts trivially on $\bar{V}$.
Since the action is defined using the projection
$\pi_1\colon\HA(\bar{V})\times\tilde G\to \HA(\bar{V})$, we have
$\pi_1(\gamma)=\id$. On the other hand,
$\pi_2(\gamma)\in\bigcap_{g\in
G}\tilde\rho(g)\tilde\Gamma\tilde\rho(g)^{-1}$. Since $(\tilde
G,\tilde\Gamma)$ is reduced, we get $\pi_2(\gamma)=e$. Thus
$\gamma=e$, and assumption (ii) of Lemma~\ref{suff_reduced} is also
satisfied.
\end{proof}

\subsection{Induced representations of $\redheck{V\rtimes G}{V_0\rtimes
\Gamma}$}\label{subsection_corner_rep} Suppose that $(V\rtimes G,
V_0\rtimes \Gamma)$ is a Hecke pair.  We denote by $\rho$ the dense
embedding of  $(V\rtimes G, V_0\rtimes \Gamma)$ in a completion of
the form $(\bar{V}\rtimes\bar{G}, \bar V_0\rtimes\bar{\Gamma})$,
which exists by Proposition~\ref{sdp_completion_universal}.

Choose left Haar measures $\mu_{\bar{V}}$  and $\mu_{\bar{G}}$ on
$\bar{V}$ and $\bar{G}$ normalized by
\begin{equation}
\mu_{\bar{V}}(\bar V_0)=1\ \text{ and }\ \mu_{\bar{G}}
(\bar{\Gamma})=1, \label{choice_of_Haar_measures}
\end{equation}
and for each $g\in \bar{G}$, let $\delta(g)$ be defined by the
formula $\mu_{\bar{V}}(\alpha_{g^{-1}} (\bar V_0))=
\delta(g)\mu_{\bar{V}}(\bar V_0)$, see for example
\cite[(15.29)]{hew-ros}.
Thus, a left Haar measure and the modular function for $\bar{V}
\rtimes \bar{G}$ are given by
$$
d\mu_{\bar{V}\rtimes \bar{G}}(vg):=\delta(g)
d\mu_{\bar{V}}(v)d\mu_{\bar{G}}(g)
$$
and
$$
\Delta_{\bar{V} \rtimes \bar{G}}(vg)= \delta(g)\Delta_{\bar{V}}(v)
\Delta_{\bar{G}}(g),
$$
for $v\in \bar{V}$ and $g\in \bar{G}$. Recall also that
$\Delta_{\bar{V} \rtimes \bar{G}} \circ \rho$ is the modular
function $\Delta_{V_0\rtimes{\Gamma}}$ associated to the given Hecke
pair. Likewise for $\Delta_{\bar{V}}$ and $\Delta_{\bar{G}}.$

Suppose now that $\chi\colon V\to \mathbb{T}$ is a character such
that $\chi(V_0)=\{1\}$. We can extend $\chi$ to a unique continuous
character on $\bar{V}$, which we continue to denote by $\chi$. The
induced representation $\operatorname{Ind}_{\bar{V}}^{\bar{V}\rtimes
\bar{G}}\chi$ acts, by definition, on the Hilbert space of functions
$\tilde{\xi}\colon \bar{V}\rtimes \bar{G}\to \mathbb{C}$ such that
$$
\tilde{\xi}(xv)=\overline{\chi(v)}\tilde{\xi}(x)\ \text{for}\
v\in\bar{V}\ \text{and}\ x\in \bar{V}\rtimes \bar{G},
$$
and $|\tilde\xi|\in L^2((\bar{V}\rtimes \bar{G})/\bar{V})$. The
representation  is simply given by left translations on this space
$$
\left(\left(\operatorname{Ind}_{\bar{V}}^{\bar{V}\rtimes
\bar{G}}\chi\right)(x) \tilde{\xi}\right)(y)=\tilde{\xi}(x^{-1}y)
\text{ for }y\in \bar{V}\rtimes \bar{G}.
$$
The map $L^2(\bar{G})\ni\xi\mapsto \tilde{\xi}$ defined by
$\tilde{\xi}(vg):=\xi(g) \overline{\chi (g^{-1}v)}$ is an
isomorphism of $L^2(\bar{G})$ onto the space of such functions,
transforming $\operatorname{Ind}_{\bar{V}}^{\bar{V}\rtimes
\bar{G}}\chi$ into the representation $\pi\colon \bar{V}\rtimes
\bar{G} \to B(L^2(\bar{G}))$ given by
\begin{equation} \label{pi_on_V_times_G}
\pi(vg)\xi(h)=\chi(\alpha_{h^{-1}}(v)) \xi(g^{-1}h)
\end{equation}
for $v\in \bar{V}$ and $g, h\in \bar{G}$.

\begin{lemma}\label{pi(p)_as_l2}
Let $\chi$ be a character of $V$ whose kernel contains $V_0$, and
let $S$ be the subset of $G$ defined by
\begin{equation} \label{eDefS}
S:=\{g\in G\mid (\chi\circ \alpha_{g^{-1}})\vert_{V_0}\equiv 1\}.
\end{equation}
Then the space of $\bar V_0\rtimes \bar{\Gamma}$-invariant vectors
for the representation $\pi$ defined by \eqref{pi_on_V_times_G}
coincides with the space of $\bar{\Gamma}$-invariant functions $f\in
L^2(\bar{G})$ with support in the closure $\bar{S}$ of $\rho(S)$ in
$\bar G$.
\end{lemma}

\begin{proof}
Since $\Gamma$ leaves $V_0$ invariant, we have $\Gamma S=S$, and
since $\bar{\Gamma}$ is an open subgroup of $\bar{G}$, this implies
that $\bar{S}$  is $\bar{\Gamma}\rho(S)$. But
$\bar{G}=\bar{\Gamma}\rho(G)$, and thus
$$
\bar{S}=\{g\in \bar{G}\mid (\chi\circ \alpha_{g^{-1}})\vert_{\bar
V_0}\equiv 1\}.
$$
Denote by  $\pi_*$ the integrated form of the representation $\pi$
from \eqref{pi_on_V_times_G}, and let $p \in C_c(\bar{V}\rtimes
\bar{G})$ be the characteristic function of the compact open
subgroup $\bar V_0\rtimes \bar{\Gamma}$, which is a self-adjoint
projection because of our normalization
\eqref{choice_of_Haar_measures}; then
$$
({\pi_*}(p)\xi)(h)= \int_{\bar V_0}\chi(\alpha_{h^{-1}}(v))
d\mu_{\bar{V}}(v) \int_{\bar{\Gamma}}\xi(g^{-1}h)d\mu_{\bar{G}}(g)
\text{ for } \xi \in L^2(\bar{G}).
$$

Since the first factor is zero or one, depending on whether $h$
belongs to $\bar{S}$ or not, we get the result.
\end{proof}

We thus obtain a $*$-representation,  which we continue to denote
by~$\pi_*$, of $pC_c(\bar{V}\rtimes \bar{G})p$ on the subspace of
$\bar \Gamma$-invariant functions in $L^2(\bar G)$ with support in
$\bar S$. By composing $\pi_*$ with the isomorphism $\Psi\colon
\hecke{V\rtimes G}{V_0\rtimes \Gamma}\to pC_c(\bar{V}\rtimes
\bar{G})p$ of Lemma~\ref{Hecke_alg_as_*_corner} we get a
representation of the Hecke algebra. A computation similar to
(\ref{how_to_get_lrr}) yields
$$
((\pi_*\circ\Psi)(f)\xi)(h)=\sum_{y=vg\in V\rtimes G/V_0\rtimes
\Gamma} \Delta_{V_0\rtimes \Gamma}(y)^{-1/2}
f(y)\chi(\alpha_{h^{-1}}(\rho(v))) \xi(\rho(g)^{-1}h)
$$
for $f\in\hecke{V\rtimes G}{V_0\rtimes \Gamma}$, $\xi\in L^2(\bar
G)$ a $\bar \Gamma$-invariant function with support in $\bar S$, and
$h\in \bar{S}$. The unitary isomorphism $(U^*
\xi)(s):=\Delta_{\bar{G}}^{1/2}(\rho(s))\xi(\rho(s))$ from the
subspace of $\bar \Gamma$-invariant functions in~$L^2(\bar G)$ with
support in $\bar S$ onto $\ell^2(\Gamma \backslash S)$ conjugates
$\pi_*\circ\Psi$ into a representation $\pi_\chi$ of
$\hecke{V\rtimes G}{V_0\rtimes \Gamma}$ on~$\ell^2(\Gamma \backslash
S)$. Since
$$
\Delta_{V_0\rtimes \Gamma}(vg)\Delta_{\Gamma}(g)^{-1}=\delta(g)
\Delta_{V_0}(v)
$$
for $vg\in V\rtimes G$, and since $(V_0\rtimes \Gamma)gv=\Gamma g
\alpha_{g^{-1}}(V_0)v$, we get that
\begin{eqnarray*}
(\pi_\chi(f)\xi)(t)&=&\sum_{vg\in
V\rtimes G /V_0\rtimes \Gamma}\delta(g)^{-1/2}\Delta_{V_0}(v)^{-1/2}
\chi(\alpha_{t^{-1}}(v))f(vg)\xi(g^{-1}t)\\
&=&\sum_{gv\in V_0\rtimes \Gamma\backslash V\rtimes
G}\delta(g)^{1/2}\Delta_{V_0}(v)^{1/2}
\overline{\chi(\alpha_{t^{-1}}(v))}f(v^{-1}g^{-1})\xi(gt)\\
&=&\sum_{g\in\Gamma\backslash G}\sum_{v\in
\alpha_{g^{-1}}(V_0)\backslash
V}\delta(g)^{1/2}\Delta_{V_0}(v)^{1/2}
\overline{\chi(\alpha_{t^{-1}}(v))}f(g^{-1}\alpha_g(v)^{-1})\xi(gt)\\
&=&\sum_{g\in\Gamma\backslash G}\sum_{v\in V_0\backslash
V}\delta(g)^{1/2}\Delta_{V_0}(v)^{1/2}
\overline{\chi(\alpha_{(gt)^{-1}}(v))}f(g^{-1}v^{-1})\xi(gt)\\
&=&\sum_{g\in\Gamma\backslash G}\sum_{v\in V_0\backslash
V}\delta(gt^{-1})^{1/2}\Delta_{V_0}(v)^{1/2}
\overline{\chi(\alpha_{g^{-1}}(v))}f(tg^{-1}v^{-1})\xi(g),
\end{eqnarray*}
for $t\in S$. Since $\xi\in\ell^2(\Gamma\backslash S)$, the last
summation is actually over $g\in \Gamma\backslash S$. We have thus
established the first claim of the next proposition.

\begin{proposition} \label{chi_extends_to_psi}
Let $\chi$ be a character of $V$ whose kernel contains $V_0$, and
let $S\subset G$ be the set defined by (\ref{eDefS}); then the
formula
$$
(\pi_\chi(f)\xi)(t)=\sum_{s\in\Gamma\backslash S}\sum_{v\in
V_0\backslash V}\delta(ts^{-1})^{-1/2}\Delta_{V_0}(v)^{1/2}
\overline{\chi(\alpha_{s^{-1}}(v))}f(ts^{-1}v^{-1})\xi(s)
$$
defines a $*$-representation $\pi_\chi\colon\hecke{V\rtimes
G}{V_0\rtimes \Gamma}\to B(\ell^2(\Gamma\backslash S))$.

If $\chi$ is weakly contained in the regular representation of
$\bar{V}$, then $\pi_\chi$ extends to a representation of
$\redheck{V\rtimes G}{V_0\rtimes \Gamma}$ on
$\ell^2(\Gamma\backslash S)$.
\end{proposition}

\bp The second assertion holds since if $\chi$ is weakly contained
in the regular representation of $\bar{V}$ then $\pi_*$ descends to
a representation of $\red{\bar{V}\rtimes\bar{G}}$. \ep

\begin{lemma}
Under the assumptions of Proposition~\ref{chi_extends_to_psi},
suppose that $V_0\subset\alpha_s(V_0)$ for every $s\in S$. Then
$S$ is a semigroup and
$$
\pi_\chi([s]_{V_0\rtimes
\Gamma})=\delta(s)^{-1/2}\lambda([s]_\Gamma)\ \text{for}\ s\in S,
$$
where $\lambda\colon\hh(S, \Gamma)\to B(\ell^2(\Gamma\backslash S))$
is the regular representation of $\hh(S, \Gamma)$ given 
by~\eqref{def_regularr_semiHecke}.
\end{lemma}

\bp Since $\ker \chi$ contains $V_0$, the identity of $G$ is in $S$,
and since $ (s_1s_2)\inv (V_0) = s_2\inv( s_1\inv(V_0)) \subset V_0
$ for $s_1, s_2 \in S$, the set $S$ is multiplicatively closed.

By \proref{smallHecke} the map $[s]_\Gamma\mapsto[s]_{V_0\rtimes
\Gamma}$ defines a homomorphism $\hecke{S}{\Gamma}\to\hecke{V\rtimes
G}{V_0\rtimes \Gamma}$. By composing $\pi_\chi$ with it we get a
representation of $\hecke{S}{\Gamma}$ on $\ell^2(\Gamma\backslash
S)$. To see that this representation matches with $\lambda$ as 
claimed, we let $s, t\in  S$, $\xi\in \ell^2(\Gamma\backslash S)$. 
Then
$$
(\pi_\chi([s]_{V_0\rtimes \Gamma})\xi)(t)=\sum_{y\in\Gamma\backslash
S}\sum_{v\in V_0\backslash
V}\delta(ty^{-1})^{-1/2}\Delta_{V_0}(v)^{1/2}
\overline{\chi(\alpha_{y^{-1}}(v))} [s]_{V_0\rtimes
\Gamma}(ty^{-1}v^{-1})\xi(y).
$$
The value of $[s]_{V_0\rtimes \Gamma}$ at $ty^{-1}v^{-1}$ is zero
unless $ty^{-1}v^{-1}\in V_0\Gamma sV_0\Gamma=\Gamma s\Gamma V_0$,
in which case $ty^{-1}\in\Gamma s\Gamma$ and $v\in V_0$. Since then
$\delta(ty^{-1})=\delta(s)$ and $\chi(\alpha_{y^{-1}}(v))=1$, we can
compute further that
$$
(\pi_\chi([s]_{V_0\rtimes \Gamma})\xi)(t)=\sum_{y\in\Gamma\backslash
S}\delta(s)^{-1/2}[s]_\Gamma(ty^{-1})\xi(y).
$$
But the last expression is  precisely
$\delta(s)^{-1/2}(\lambda([s]_\Gamma)\xi)(t)$, and the lemma
follows. \ep

In other words, the map
$[s]_\Gamma\mapsto\delta(s)^{1/2}[s]_{V_0\rtimes \Gamma}$ from
Remark~\ref{another_homo} gives an embedding of $\hecke{S}{\Gamma}$
into $\hecke{V\rtimes G}{V_0\rtimes \Gamma}$ such that the
restriction of the representation $\pi_\chi$ to $\hecke{S}{\Gamma}$
coincides with the regular representation $\lambda$ of $\hh(S,
\Gamma).$

\medskip

Our next immediate goal is to apply the preceding considerations to
the Hecke pair $(P,P_0)$ obtained by setting $V=\mq$, $V_0=\mz$,
$G=\glq$ and $\Gamma=\slz$.  Since the Hecke pairs $(\mq, \mz)$ and
$(\glq, \slz) $ are unimodular, it follows from
Proposition~\ref{HeckeSmallinsideBig} and
Proposition~\ref{MainHecke} that
$$
\delta(g)=\Delta_{P_0}(g)=\det(g)^{-2}\text{ for }g\in G.
$$
We now fix a character $\chi$ on $\af$ such that the corresponding
pairing $ (x, y) \mapsto \chi^y(x): = \chi (xy)$ for $x,y\in \af$
gives a self-duality isomorphism $y \mapsto \chi^y $ of $\af$ to
$\widehat{\af}$ in which  $\rr$ corresponds to $\rr^\perp $. Thus in
particular $\chi(\rr)=\{1\}$. Denote by $\Tr $ the usual trace on
$2\times 2$ matrices; then a similar pairing  given by
\begin{equation}
(a,b) \mapsto \chi^b(a):=\chi(\Tr (ab)),\  a,b \in \ma,
\label{def_of_chi_on_matrices_on_af}
\end{equation}
implements  a self-duality isomorphism $a \mapsto \chi^a$ of  $\ma$
onto $\widehat\ma$  in which $\mr$ corresponds to $\mr^\perp$. This
pairing is noncanonical, but for any choice of $\chi$ as above we
have that
\[
\chi^n(\alpha_g(m)) = \chi(\Tr(mg^{-1}n)) = \chi^{g^{-1}n}(m)
\]
for $g\in \glap$, and this says that the action by right
multiplication by $g$ on $\ma$ is transformed by the self-duality
into left multiplication.

An element $g\in\glq$ belongs to
$S:=\{h\in\glq\mid(\chi\circ\alpha_{h^{-1}})|_{\mz}\equiv1\}$
precisely when $gw\in\mr^\perp=\mr$ for all $w\in\mr$, that is, when
$g\in\glq\cap\mr=\mzp$. We thus get a representation of our Hecke
algebra on the space $\ell^2(\Gamma\backslash S)$. This
representation is faithful by the following known result.

\begin{lemma}
Let $G$ be a locally compact group acting on a locally compact space
$X$, and suppose that $x$ is a point in $X$ with dense orbit. Let
$\chi_x$ denote evaluation at $x$ on $C_0(X)$. Then
$\operatorname{Ind}_{C_0(X)}^{C_0(X)\rtimes_r G}\chi_x$ is faithful.
\label{lemma_psi_faithful}
\end{lemma}

\bp Denote by $\pi$ the representation $\oplus_{g\in G}\chi_{gx}$ of
$C_0(X)$. The representations
$\operatorname{Ind}_{C_0(X)}^{C_0(X)\rtimes_r G}\chi_{gx}$ and
$\operatorname{Ind}_{C_0(X)}^{C_0(X)\rtimes_r G}\chi_x$ are
equivalent for every $g$ in $G$, and therefore
$\operatorname{Ind}_{C_0(X)}^{C_0(X)\rtimes_r G}\chi_x$ is
quasi-equivalent to the representation
$\operatorname{Ind}_{C_0(X)}^{C_0(X)\rtimes_r G}\pi$, which is
faithful because $\pi$ is faithful. The lemma follows. \ep

We want to apply Lemma~\ref{lemma_psi_faithful} to the group $G =
\glap$ acting on $\hat V=\ma$ by multiplication on the left. The
density of $\mzp$ in $\mtwo(\rr)$ implies that $\glap w$ is dense in
$\ma$ for every $w\in\glr$. Therefore the lemma applies, and we get
the following.

\begin{theorem} \label{pi_faithful_our_pair}
For the Hecke pair $(P,P_0)=(\mq\rtimes \glq, \mz\rtimes \slz)$, the
semigroup $S=\mzp$, and each $w\in\glr$,  the representation
$\pi_w:=\pi_{\chi^w}$ constructed as in
Proposition~\ref{chi_extends_to_psi} from the character $\chi^w$
defined by \eqref{def_of_chi_on_matrices_on_af} is a faithful
representation of $\redheck{P}{P_0}$ on $\ell^2(\Gamma \backslash
S)$. The restriction of this representation to the subalgebra
$\hh(\mzp,\slz)$ is simply the left regular representation rescaled
by a factor given by the determinant:
$$
\pi_w([s]_{P_0})=\det(s)\lambda([s]_\Gamma) \text{ for }s\in \mzp.
$$
\end{theorem}

\bigskip

\section{The symmetry group and the symmetric part of the algebra}
\label{symmetries}

In this section we continue our study of the reduced Hecke
C$^*$-algebra of the pair
$$
(P,P_0) = \left(\matr{\mq}{\glq}, \matr{\mz}{\slz} \right),
$$
with the aim of showing that there exists an action of the compact
group $\rr^*$ of invertible elements in~$\rr$  such that the fixed
point algebra decomposes into a tensor product of algebras
corresponding to different primes.

\subsection{The crossed product picture of $\redheck{P}{P_0}$ and symmetries}

Note that if we have a Hecke pair $(N,N_0)$, then every automorphism
$\theta$ of $N$ leaving~$N_0$ invariant defines an automorphism
$\alpha_\theta$ of the Hecke algebra by
$\alpha_\theta(f)=f\circ\theta^{-1}$. Moreover, since any completion
of $(N,N_0)$ defines an isomorphic Hecke algebra, we can consider
automorphisms of completions. As  showed in \secref{scompletion},
the Schlichting completion of our pair $(P,P_0)$ is
$$
(\bar P, \bar P_0) = \left(\matr{\ma}{\glap},
\matr{\mr}{\slr}\right).
$$
For $r\in\G(\rr)$ the automorphism of $\bar P$ defined by
$$
\bar{P}\ni x\mapsto\matr{0}{r}x\matr{0}{r^{-1}}
$$
leaves $\bar P_0$ invariant, and hence defines an automorphism
$\alpha_r$ of $\hh(\bar{P}, \bar{P_0})$.  Since $\hh(\bar{P},
\bar{P_0})$ consists by definition of $\bar P_0$-biinvariant
functions, $\alpha_r$ is trivial for $r\in\slr$. Thus we get an
action of the group $\G(\rr)/\slr\cong\rr^*$ on $\hh(\bar{P},
\bar{P_0})$ and on $\hh({P}, {P_0})$.

Recall from \secref{stuff_compl} that the reduced Hecke
C$^*$-algebra $\redheck{P}{P_0}$ is isomorphic to the corner in
$\red {\bar P}$ defined by the projection corresponding to the
compact open subgroup $\bar P_0$. It is standard that the
C$^*$-algebra of a semidirect product can be written as a crossed
product, so $ \red{\bar P} \cong \red{\ma} \rtimes \glap$, where the
action of $\glap $ on $\red{\ma} $ is obtained from right
multiplication of $\glap $ on $\ma$. But the self-duality of
$\mtwo(\af)$ discussed before Lemma~\ref{lemma_psi_faithful} allows
us to transpose the action of $\glap$ into left multiplication  on
$\ma$. Hence using the Fourier transform and the transposed action
we may write
\[
\red{\bar P}\cong \red{\ma} \rtimes \glap \cong C_0(\ma) \rtimes_r
\glap
\]
where the action of $\glap$ on $C_0(\ma)$ is given by $(g\cdot
f)(m)=f(g^{-1}m)$ for $g\in \glap$ and $m\in \ma$. Moreover, since
$\mtwo(\rr)^\perp=\mtwo(\rr)$, this isomorphism carries the
projection corresponding to $\bar P_0$ into
$$
p_0=\31_{\mtwo(\rr)}\int_{\slr}\lambda_gdg,
$$
where $\31_{\mtwo(\rr)}$ is the characteristic function of the set
$\mtwo(\rr)$, and $\lambda_g$ the element corresponding to
$g\in\slr$ in the multiplier algebra of the second crossed product.
Thus, with $\Psi$ from  \lemref{Hecke_alg_as_*_corner}   we have
\begin{equation} \label{eiso}
\Psi:\redheck{P}{P_0}\overset{\cong}{\longrightarrow}
p_0(C_0(\ma)\rtimes_r\glap)p_0.
\end{equation}

We can regard $C_0(\ma)\rtimes_r\glap$ as a completion of the
algebra of continuous complex-valued functions with  compact support
on $\glap\times\ma$, endowed with the convolution product
\begin{equation}
(f_1*f_2)(g,m)=\int_{\glap}f_1(gh^{-1},hm)f_2(h,m)dh
\label{conv_prod_bigcp}
\end{equation}
and involution
\begin{equation}
f^*(g,m)=f(g^{-1},gm).
\label{inv_bigcp}
\end{equation}
Cutting down to the corner determined by $p_0$ has two effects on
functions $f$ on the space $\glap \times \ma$. The first one is that
it reduces the support to pairs $(g,m)\in\glap\times\ma$ such that
both $m$ and $gm$ are in $\mr$ and second, it forces the invariance
under the action $\slr\times\slr$ on $\glap\times\ma$ given by
\begin{equation}
(g_1,g_2)(g,m)=(g_1gg_2^{-1},g_2m)\text{ for }g_1, g_2\in \slr.
\label{action_sl2r_x_sl2r}
\end{equation}

Therefore, if we denote by $\fibb{\glap}{\slr}{\mr}$ the quotient of
the space
$$
\{(g,m)\in\glap\times\mtwo(\rr)\,\mid\, gm\in\mtwo(\rr)\}
$$
by the above  action  of $\slr\times\slr$, then
$p_0(C_0(\ma)\rtimes_r\glap)p_0$ is a completion of the algebra
$C_c(\fibb{\glap}{\slr}{\mr})$ of compactly supported  continuous
functions on $\fibb{\glap}{\slr}{\mr}$ with convolution and
involution given by \eqref{conv_prod_bigcp} and~\eqref{inv_bigcp}.
Since $\slr\backslash\glap=\Gamma\backslash\glq$ by
Lemma~\ref{glqdenseinglap}, we can refine formula
\eqref{conv_prod_bigcp} for the convolution of two such functions to
\begin{equation}\label{new_conv_prod}
(f_1*f_2)(g,m)
=\sum_{h\in\Gamma\backslash\glq:\,hm\in\mr}f_1(gh^{-1},hm)f_2(h,m).
\end{equation}
If we regard elements of $\hecke{P}{P_0}$ as functions defined on
$\glq\times\mq$ (it is convenient now to use this order), the
isomorphism $\Psi$ from  \eqref{eiso} obtained by applying the
Fourier transform in the second variable and multiplying by
$\Delta_{P_0}^{-1/2}$ (due to Lemma~\ref{Hecke_alg_as_*_corner}) is
then explicitly given by
\begin{equation} \label{FT_f}
\Psi(f)(g,m)=\det(g)\sum_{n\in\mq/\mz} f(g,n)\chi(\Tr(nm)).
\end{equation}
For the rest of the paper we shall identify $\redheck{P}{P_0}$ with
its image under $\Psi$. So with this convention, the automorphism
$\alpha_r$ is now given by
\begin{equation}\label{Ft_alpha}
\alpha_r(f)(g,m)=f(r^{-1}gr,r^{-1}m),
\end{equation}
for $f\in C_c(\fibb{\glap}{\slr}{\mr})$.

\subsection{Tensor product decomposition of the fixed point algebra}
To understand the fixed point algebra of the action
\eqref{Ft_alpha}, for each prime $p$  define a subgroup $G_p$ of
$\glq$ by
$$
G_p:=\glp(\Z[p^{-1}])=\G(\Z[p^{-1}])\cap\glq.
$$
It is not difficult to show, see \cite[Section~3]{LLN}, that $G_p$
is the subgroup of $\glq$ generated by $\Gamma$ and $\matr{0}{p}$,
and also that it is the group of elements $g\in\glq$ such that the
image of $g$ in $\G(\Q_q)$ lies in $\G(\Z_q)$ for $q\ne p$. Consider
now the Hecke pair
\begin{equation}\label{p_component_of_our_hecke}
(\mtwo(\Z[p^{-1}])\rtimes G_p,\mz\rtimes\Gamma)=
\left(\matr{\mtwo(\Z[p^{-1}])}{\glp(\Z[p^{-1}])}, \matr{\mz}{\slz}
\right).
\end{equation}
Similarly to Proposition~\ref{completion_our_pair}, the Schlichting
completion of \eqref{p_component_of_our_hecke} is the pair
$$
(\mtwo(\Q_p)\rtimes\glp(\Q_p),\mtwo(\Z_p)\rtimes\sltwo(\Z_p)),
$$
where $\glp(\Q_p)$ is the group of elements $g\in\G(\Q_p)$ such that
$\det(g)$ is an integral power of~$p$. The reduced Hecke
C$^*$-algebra of this pair can be regarded as a subalgebra of
$\redheck{P}{P_0}$. If we consider $\redheck{P}{P_0}$ as a
completion of $C_c(\fibb{\glap}{\slr}{\mr})$ as explained above,
then this subalgebra has the following description.

\begin{lemma}\label{gothic_Ap}
The reduced Hecke C$^*$-algebra
$\mathfrak{A}_p:=\redheck{\mtwo(\Z[p^{-1}])\rtimes
G_p}{\mz\rtimes\Gamma}$ of the pair~\eqref{p_component_of_our_hecke}
is the closure of the space of $(\slr\times\slr)$-invariant
compactly supported continuous functions~$f$ on
$$
\{(g,m)\in\glap\times\mtwo(\rr)\,\mid\, gm\in\mtwo(\rr)\}
$$
such that $f(g,m)=0$ if $g\notin G_p\slr$, and $f(g,m)=f(g,m')$ if
$m_p=m'_p$.
\end{lemma}

\begin{proof}
The closure of $G_p$ in $\glap$ is the group $G_p\slr$, and the
closure of $\mtwo(\Z[p^{-1}])$ in $\ma$ is
$\mtwo(\Q_p\times\prod_{q\ne p}\Z_q)$. The annihilator of
$\mtwo(\Q_p\times\prod_{q\ne p}\Z_q)$ is the group of elements
$m\in\mtwo(\af)$ such that $m_p=0$ and $m_q\in\mtwo(\Z_q)$ for $q\ne
p$. Since a function on $\ma/\mr$ is supported on
$\mtwo(\Q_p\times\prod_{q\ne p}\Z_q)$ if and only if its Fourier
transform is invariant under the translations by elements of the
annihilator, and the latter means that the value of the Fourier
transform at $m\in\mr$ depends only on $m_p$, the result follows.
\end{proof}

Note that the action of $\alpha_r$ on $\mathfrak{A}_p$ depends only
on $r_p$. Indeed, assume $r_p=1$ in $\gl(\Z_p)$. Since every double
coset of $\Gamma$ in $\glq$ has a diagonal representative, see e.g.
\cite[Chapter IV]{kri}, it suffices to compute $\alpha_r(f)(g,m)$
for $g\in G_p$ diagonal. Since $\alpha_r$ depends only on $\det(r)$,
we may assume that $r$ is also diagonal, and since
$(r^{-1}m)_p=m_p$, we get
$$
\alpha_r(f)(g,m)=f(g,r^{-1}m)=f(g,m).
$$
Therefore the action of $\rr^*$, when restricted to
$\mathfrak{A}_p$, defines an action of $\Z_p^*$. Alternatively, this
action can also be obtained from conjugation of elements in
$\matr{\mtwo(\Q_p)}{\glp(\Q_p)}$ by matrices $\diag{1}{r}$ with
$r\in\G(\Z_p)$.

\smallskip

We can now formulate the main result of the section.

\begin{theorem} \label{Symm}
The subalgebras $\mathfrak{A}_p^{\Z^*_p}={\redheck{\mtwo(\Z[p^{-1}])
\rtimes G_p}{\mz\rtimes\Gamma}}^{\Z^*_p}$ of $\redheck{P}{P_0}$
mutually commute for different primes, and
$$
{\redheck{P}{P_0}}^{\rr^*}\cong\bigotimes_{p\in\primes}
\mathfrak{A}_p^{\Z^*_p}.
$$
\end{theorem}

In the proof of this theorem  we rely on lattices and their
properties, and since we will also use these later we collect the
data we need in the next remark.

\begin{remark}\label{gen_rem_lattices}
A lattice $L$ in $\R^2$ is commensurable with $\Z^2$ if it is
contained in $\Q^2$. We identify $\Gamma\backslash S$, where
$S=\mzp$, with the set of lattices in $\R^2$ containing $\Z^2$:
namely, for $s\in S$ we let $L=s^{-1}\Z^2$. Equivalently, we can
consider lattices in $\af^2$ containing $\rr^2$ as follows: given a
lattice $L\subset\Q^2$, the closure $\bar L$ of $L$ in $\af^2$ is a
lattice, and by \cite[Theorem V.2]{weil}, the map $L\mapsto\bar L$
is a bijection between lattices in $\Q^2$ and lattices in $\af^2$
with inverse given by $\bar L\mapsto \Q^2\cap \bar L$. The group
$\G(\rr)$ acts on the space of lattices by
$$
rL:=r\bar L\cap\Q^2\text{ for }r\in\G(\rr).
$$

Given a lattice $L$ commensurable with $\Z^2$, we denote by $L_p$
the lattice $\Q^2\cap(\Z_p\otimes_\Z L)$: in other words, $L_p$ is
the unique lattice in $\Q^2$ such that the closure of $L_p$ in
$\Q_p^2$ coincides with that of $L$, and the closure of $L_p$ in
$\Q_q^2$ is $\Z_q^2$ for $q\ne p$. If $L$ contains $\Z^2$ and we
write $L_p=s^{-1}\Z^2$ and use that the closure of $L_p$ in $\Q_q^2$
is $\Z_q^2$ for $q\ne p$, we deduce that $s\in \G(\Z_q)$ for $q\ne
p$. Hence $s\in S_p:=G_p\cap S$. Note that $S_p$ is the semigroup of
matrices $m\in\mzp$ such that $\det(m)$ is a nonnegative power
of~$p$. Thus to every $L\in \Gamma \backslash S$ corresponds a
family $(L_p)_{p\in \mathcal{P}}$ with $L_p\in \Gamma \backslash
S_p$. An equivalent description of $\Gamma\backslash S_p$ is as the
set of lattices $L$ such that $\Z^2\subset L$ and $L/\Z^2$ is a
$p$-group, i.e. its order is a power of~$p$. The map
$\delta_L\mapsto\otimes_{p\in\primes}\delta_{L_p}$ defines a unitary
isomorphism
\begin{equation}\label{tensor_decomp_spaces}
(\ell^2(\Gamma\backslash
S),\delta_{\Z^2})\cong\otimes_{p\in\primes}(\ell^2(\Gamma\backslash
S_p),\delta_{\Z^2}).
\end{equation}

Finally, since $\slr\backslash\glap=\Gamma\backslash\glq$ can be
identified with the set of lattices commensurable with $\Z^2$, a
function on $\fibb{\glap}{\slr}{\mr}$ can be thought of as a
function~$f$ on the set of pairs $(L,m)$, where $L$ is a lattice
commensurable with $\Z^2$ and $m\in\mr$ is a matrix whose columns
belong to $\bar L\subset\af^2$, such that $f$ is invariant under the
action of $\slr$ given by $\gamma(L,m)=(\gamma L,\gamma m)$. In this
picture the action of $\G(\rr)$ on $\redheck{P}{P_0}$ is given by
$\alpha_r(f)(L,m)=f(r^{-1}L,r^{-1}m)$.
\end{remark}

\bp[Proof of \thmref{Symm}] Denote by $\pi$ the representation of
$\redheck{P}{P_0}$ on $\ell^2(\Gamma\backslash S)$ defined by $w=1$
as described in Theorem~\ref{pi_faithful_our_pair}. Using the
identification of the reduced Hecke C$^*$-algebra with the corner of
the crossed product via \eqref{eiso}, for $f\in
C_c(\fibb{\glap}{\slr}{\mr})$ and $s\in S$ we get
$$
\pi(f)\delta_{\Gamma s}=\sum_{t\in\Gamma\backslash S}
f(ts^{-1},s)\delta_{\Gamma t}.
$$
If we view $f$ as a function on a pair $(L,m)$ as explained in
Remark~\ref{gen_rem_lattices}, for $L=s^{-1}\Z^2$ we obtain
\begin{equation} \label{esymm1}
\pi(f)\delta_L=\sum_{L'\supset\Z^2}f(sL',s)\delta_{L'}.
\end{equation}

For each prime $p$ we have a similar representation $\pi_p$ of
$\mathfrak{A}_p$ on $\ell^2(\Gamma\backslash S_p)$.  Let $f\in
\mathfrak{A}_p$ and view it, by Lemma~\ref{gothic_Ap} and
Remark~\ref{gen_rem_lattices}, as an $\slr$-invariant function on
$(L,m)$ such that $f(L,m)$ is nonzero only if $L$ is defined by an
element in $G_p$, and $f(L,m)=f(L,m')$ when $m_p=m'_p$. Then for
$L''=t^{-1}\Z^2$ with $t\in S_p$, we have an analogue of
\eqref{esymm1} for $\pi_p$:
\begin{equation} \label{esymm2}
\pi_p(f)\delta_{L''}=\sum_{ \substack{L'\supset\Z^2:\\L'/\Z^2\
\text{ is a} \ p\text{-group}}}f(tL',t)\delta_{L'}.
\end{equation}

We claim that for every $f\in{\mathfrak{A}_p}^{\Z_p^*}$ we have
$\pi(f)=\pi_p(f)\otimes1$ with respect to the decomposition of
$\ell^2(\Gamma \backslash S)$ into the tensor product of
$\ell^2(\Gamma\backslash S_p)$ and $\otimes_{q\neq p}(\ell^2(\Gamma
\backslash S_q), \delta_{\Z^2})$ given by
\eqref{tensor_decomp_spaces}. Indeed, fix such $f$, and let
$L=s^{-1}\Z^2$ for $s\in S$ be a lattice containing $\Z^2$. In the
right hand side of  \eqref{esymm1}, the value $f(sL',s)$ is nonzero
only if $sL'$ is defined by an element in $S_p$, that is,
$(sL')_q=\Z^2$ for $q\ne p$, or equivalently,
$L'_q=(s^{-1}\Z^2)_q=L_q$. Then the summation is over lattices $L'$
for which we have a decomposition $\delta_{L'}=\delta_{L'_p}\otimes
\bigotimes_{q\neq p}\delta_{L_q}$, and (\ref{esymm1}) becomes
$$
\pi(f)\delta_L =\sum_{\substack{L'\supset\Z^2:\\ L'_q=L_q\
\text{for} \ q\ne p}}f(sL',s)\delta_{L'}
=\left(\sum_{\substack{L'\supset\Z^2:\\ L'_q=L_q\ \text{for} \ q\ne
p}}f(sL',s)\delta_{L'_p}\right)\otimes \bigotimes_{q\neq
p}\delta_{L_q}.
$$
Comparing the last sum with (\ref{esymm2}) (with $L''=L_p$) and
using that when $L'$ runs through all lattices containing $\Z^2$
such that $L'_q=L_q$ for $q\ne p$ then $L'_p$ runs through all
lattices such that $L'_p/\Z^2$ is a $p$-group, we see that to prove
the claim it suffices to check that
\begin{equation}\label{on_the_way}
f(sL',s)=f(tL'_p,t)\ \ \text{if}\ \ L=s^{-1}\Z^2,\ L_p=t^{-1}\Z^2\
\text{and}\ L'_q=L_q\ \text{for}\ q\ne p.
\end{equation}
Since $f$ is $\G(\rr)$-invariant and $f(sL',s)=f(sL',m)$ if $s=m_p$
in $\mtwo(\Z_p)$, equality \eqref{on_the_way} will hold if we find
an element $r\in\G(\rr)$ such that $s=r_pt$ in $\G(\Q_p)$ and
$sL'=rtL'_p$.

We take $r_p=st^{-1}$ and $r_q=1$ for $q\ne p$. Since the closures
of $L$ and $L_p$ in $\Q_p^2$ coincide by definition, we have
$s^{-1}\Z^2_p=t^{-1}\Z^2_p$, so that $r_p=st^{-1}\in\G(\Z_p)$. Since
the closures of $L'$ and $L'_p$ in~$\Q_p^2$ coincide, it is also
clear that the closures of $sL'$ and $rtL'_p$ in $\Q^2_p$ coincide.
On the other hand, for $q\ne p$ we have
$(sL')_q=\Z^2=(tL'_p)_q=(rtL'_p)_q$, so that the closures of $sL'$
and $rtL'_p$ in $\Q_q^2$ coincide. Hence $sL'=rtL'_p$, and the claim
is proved.

Since the representation $\pi$ is faithful, the claim implies that
the algebras $\mathfrak{A}_p^{\Z^*_p}$ mutually commute for
different primes, and the C$^*$-subalgebra they generate is
isomorphic to their tensor product. It remains to show that this
subalgebra coincides with the whole fixed point algebra.

Towards this end we first prove that if $f_p\in
\mathfrak{A}_p^{\Z_p^*}$ for $p$ in a finite set $F$ of prime
numbers, then
\begin{equation} \label{esymm3}
(*_{p\in F}f_p)(L,m)=\begin{cases}\prod_{p\in F}f_p(L_p,m), &
\text{if}\ L_q=\Z^2\ \text{for}\ q\notin F,\\
0, &\text{otherwise}.
\end{cases}
\end{equation}
This can be deduced from the equality $\pi(f)=\pi_p(f)\otimes 1$
above, but we can give a direct argument as follows. To simplify
notation we only consider the case of a two-point set. So assume
$F=\{p_1,p_2\}$. For $f_1$ and $f_2$ as in \eqref{esymm3} we have
$$
(f_{p_1}*f_{p_2})(L,m)=\sum_{h\in\Gamma\backslash\glq:\,hm\in\mr}
f_{p_1}(hL,hm)f_{p_2}(h^{-1}\Z^2, m).
$$
The second factor in the right hand side can be nonzero only if
$h\in G_{p_2}$. But for the first factor to be nonzero we need
$(hL)_q=\Z^2$ for $q\ne p_1$. In particular, $L_q=\Z^2$ for $q\ne
p_1,p_2$, and $L_{p_2}=h^{-1}\Z^2$, which shows that there is at
most one nonzero summand. The element $r$ defined by $r_{p_1}=h$ and
$r_q=1$ for $q\ne p_1$ lies in $\G(\rr)$, and $hL=rL_{p_1}$ because
both lattices have the same closure in $\Q_{p_1}^2$ by the choice of
$r_{p_1}$, and the same closure  $\Z_q^2$ in $\Q_q^2$ for $q\neq
{p_1}$. This implies that $f_{p_1}(hL,hm)=f_{p_1}(rL_{p_1}, rm)$,
and this is $f_{p_1}(L_{p_1}, m)$ by $\G(\rr)$-invariance, proving
(\ref{esymm3}).

Now let $L\supset\Z^2$ be a lattice, and $F\subset\primes$ a finite
subset such that $L_q=\Z^2$ for $q\notin F$. Let $U$ be an open
compact subset of $\mtwo(\rr)$ of the form
$$
\prod_{p\in F}U_p\times\prod_{q\notin F}\mtwo(\Z_q)
$$
such that the column vectors of any element $m\in U$ belong to $\bar
L$. Let $f$ be the characteristic function of the set
$\G(\rr)(L,U)=\{(rL,rm)\,\mid\,r\in\G(\rr), \ m\in U\}$. The linear
span of such functions is dense in $\redheck{P}{P_0}^{\rr^*}$. On
the other hand, denoting by $f_p$ the characteristic function of the
set $\G(\rr)(L_p,U_p\times\prod_{q\ne p}\mtwo(\Z_q))$ we have
$f_p\in\mathfrak{A}_p^{\Z_p^*}$, and $f=*_{p\in F}f_p$
by~(\ref{esymm3}). This completes the proof of the theorem. 
\ep

\begin{remark}\label{pi_w_equiv}
For  each  $w\in\G(\rr)$ consider the representation $\pi_w$ of
$\redheck{P}{P_0}$ on $\ell^2(\Gamma\backslash S)$ defined in
Theorem~\ref{pi_faithful_our_pair}. As opposed to the one
dimensional case, the restrictions of these representations to
$\redheck{P}{P_0}^{\rr^*}$ still depend on $w$. Let us show that
nevertheless, they are equivalent.

Similarly to \eqref{esymm1}, for $L=s^{-1}\Z^2$ and $f\in
C_c(\fibb{\glap}{\slr}{\mr})$ we have
\begin{equation} \label{erep}
\pi_w(f)\delta_L=\sum_{L'\supset\Z^2}f(sL',sw)\delta_{L'}.
\end{equation}

Define a unitary $U_w$ on $\ell^2(\Gamma\backslash S)$ by
$U_w\delta_L=\delta_{wL}$, and let $s_w\in S$ be such that
$w^{-1}L=s_w^{-1}\Z^2$. Then
$$
U_w\pi(f)U^*_w\delta_L=U_w\pi(f)\delta_{w^{-1}L}
=\sum_{L'\supset\Z^2}f(s_wL',s_w)\delta_{wL'}
=\sum_{L'\supset\Z^2}f(s_ww^{-1}L',s_w)\delta_{L'}
$$
Since $w^{-1}s^{-1}\Z^2=w^{-1}L=s_w^{-1}\Z^2$, we have $sw=rs_w$ for
some $r\in\G(\rr)$. Therefore for $f\in\redheck{P}{P_0}^{\rr^*}$ we
get
$$
f(s_ww^{-1}L',s_w)=f(r^{-1}sL',r^{-1}sw)=f(sL',sw).
$$
Thus $U_w\pi(f)U^*_w=\pi_w(f)$ for $f\in\redheck{P}{P_0}^{\rr^*}$.
Note that if $\det(w)=1$ then $r\in\slr$, and hence
$U_w\pi(\cdot)U_w=\pi_w$ as representations of $\redheck{P}{P_0}$.
\end{remark}

\subsection{The structure of the fixed point algebra}
The fixed point algebra $\redheck{P}{P_0}^{\rr^*}$ contains two
important subalgebras:  the diagonal subalgebra and the semigroup
Hecke algebra. We start by describing the diagonal subalgebra. Due
to the equality
$$
p_0C_0(\ma)p_0=C(\slr\backslash\mtwo(\rr))p_0,
$$
the algebra $C(\slr\backslash\mtwo(\rr))$ can be viewed as a
subalgebra of $p_0(C_0(\mtwo(\af)\rtimes\glap)p_0$. With our
equivalent picture of functions on $\fibb{\glap}{\slr}{\mr}$, we are
looking at ($\slr\times\slr$)-invariant functions on
$\slr\times\mtwo(\rr)$, which depend only on the second coordinate.
Then
$C(\G(\rr)\backslash\mtwo(\rr))=C(\slr\backslash\mtwo(\rr))^{\rr^{*}}$
is a subalgebra of $\redheck{P}{P_0}^{\rr^*}$. Likewise,
$C(\G(\Z_p)\backslash \mtwo(\Z_p))$ is a subalgebra of
$\mathfrak{A}_p^{\Z_p^*}$ for each prime $p$.

The algebra $C(\G(\rr)\backslash\mtwo(\rr))$ can be described in
terms of lattices as follows. For $L=s^{-1}\Z^2$ with $s\in S$, we
denote by $\pi_L$ the characteristic function of the set
$\mtwo(\rr)s$. Then $\pi_L$ belongs to
$C(\G(\rr)\backslash\mtwo(\rr))$.

\begin{proposition}
The algebra $C(\G(\rr)\backslash\mtwo(\rr))$ is the universal
C$^*$-algebra generated by projections $\pi_L$, $L\supset\Z^2$,
satisfying the relations $\pi_L\pi_{L'}=\pi_{L+L'}$.
\end{proposition}

\bp We first check that $\pi_L\pi_{L'}=\pi_{L+L'}$ in
$C(\G(\rr)\backslash\mtwo(\rr))$. Note that if $L=s^{-1}\Z^2$ then
$m\in\mtwo(\rr)s$ if and only if $ms^{-1}\in\mtwo(\rr)$, or
equivalently, $ms^{-1}\rr^2\subset\rr^2$, that is, $m\bar
L\subset\rr^2$. Since $\overline{L+L'} =\bar L+\bar L'$, the
identity $\pi_L\pi_{L'}=\pi_{L+L'}$ follows.

Denote by $B$ the universal C$^*$-algebra generated by projections
$\chi_L$, $L\supset\Z^2$, satisfying the relations
$\chi_L\chi_{L'}=\chi_{L+L'}$. Then we have a homomorphism $B\to
C(\G(\rr)\backslash\mtwo(\rr))$ which maps $\chi_L$ onto~$\pi_L$.
Since any finite set of the $\chi_L$'s generate a finite dimensional
algebra, to prove injectivity it is enough to check that the map is
injective on the linear span of the $\chi_L$'s. Assume
$\sum_i\lambda_i\pi_{L_i}=0$. Choose~$i_0$ such that $L_i$ is not
contained in $L_{i_0}$ for~$i\ne i_0$, and write
$L_{i_0}=s^{-1}\Z^2$ for $s\in S$. Then $\pi_{L_{i_0}}(s)=1$ and
$\pi_{L_i}(s)=0$ for~$i\ne i_0$, since $s\bar L_i$ is not contained
in $\rr^2$. Hence $\lambda_{i_0}=0$. Thus we inductively get
$\lambda_i=0$ for all~$i$.

To prove surjectivity we have to check that the functions $\pi_L$
separate points of $\G(\rr)\backslash\mtwo(\rr)$. For this we claim
that $m,n\in\mtwo(\rr)$ belong to the same $\G(\rr)$-orbit if and
only if $m^t\rr^2=n^t\rr^2$, where $m^t$ and $n^t$ denote the
transposed matrices. Indeed, assume $m^t\rr^2=n^t\rr^2$. Then for
each $p$ the $\Z_p$-module $m^t_p\Z_p^2$ is finitely generated and
torsion-free, hence it is free. Therefore $\Z_p^2$ decomposes into
the direct sum of $\Z_p^2\cap\Ker m^t_p$ and a module $V'_p$ such
that $V'_p\cap\Ker m^t_p=0$. Similarly find a module $V''_p$ for
$n$. Since $m^t_p\Z_p^2=n^t\Z_p^2$, the ranks of $V'_p$ and $V''_p$
coincide. Define $r_p\in\G(\Z_p)$ such that $r_p(\Ker
n^t_p\cap\Z_p^2)=\Ker m^t_p\cap\Z_p^2$, $r_pV''_p=V'_p$ and
$m^t_pr_p=n^t_p$. Then $m^tr=n^t$, so that $n=r^tm$. The converse
statement is clear.

Now if $m$ and $n$ belong to different $\G(\rr)$-orbits, we may
assume that $n^t\rr^2$ is not contained in~$m^t\rr^2$. Then there
exists a lattice in $\rr^2$ which contains $m^t\rr^2$ but does not
contain $n^t\rr^2$, and so the dual lattice $\bar L$ has the
properties that $\rr^2\subset \bar L$, $m\bar L\subset\rr^2$ but
$n\bar L$ is not contained in $\rr^2$. In other words, $\pi_L(m)=1$
and $\pi_L(n)=0$, and surjectivity is proved.\ep

We note that in the representation $\pi$ considered in
\eqref{esymm1}, the projections $\pi_L$ act as follows:
$$
\pi(\pi_L)\delta_{L'}=\begin{cases}
\delta_{L'},&\text{if}\ L\subset L',\\
0,&\text{otherwise}.
\end{cases}
$$

To introduce the second important subalgebra of
$\redheck{P}{P_0}^{\rr^{*}}$, recall that by \proref{smallHecke} and
the remark following it, the map
$[s]_\Gamma\mapsto\det(s)^{-1}[s]_{P_0}$ gives an embedding of the
classical Hecke algebra $\hecke{S}{\Gamma}=\hecke{\mzp}{\slz}$ into
$\hecke{P}{P_0}$. Note that the right hand side of \eqref{FT_f} with
$f=\det(s)^{-1}[s]_{P_0}$ is the characteristic function of the set
$\slr s\slr\times\mtwo(\rr)$. In terms of pairs $(L,m)$, the element
$[s]_\Gamma$ corresponds to the characteristic function of the set
$\Gamma L_s\times\mtwo(\rr)$, where $L_s=s^{-1}\Z^2$. Since $s$ can
be assumed to be diagonal, we see that  $\hecke{S}{\Gamma}$ is
contained in $\redheck{P}{P_0}^{\rr^*}$.

The algebra $\hecke{S}{\Gamma}$ is the tensor product of its
subalgebras $\hecke{S_p}{\Gamma}$, see \cite[Chapter IV]{kri}. Our
\thmref{Symm} can be regarded as a generalization of this fact.
Moreover, by e.g. \cite{kri}, $\hecke{S_p}{\Gamma}$ is generated by
two elements
$$
u_p=\left[\diag{p}{p}\right]\ \ \text{and}\ \
v_p=\left[\diag{1}{p}\right].
$$
Since $\Gamma\diag{1}{p}\Gamma$ is the set of all matrices in $\mzp$
with determinant $p$, the element $v_p$ is the characteristic
function of the set of pairs $(L,m)$ such that $|L/\Z^2|=p$. Thus
for $s\in S$ we have $v_p(sL',sm)=1$ if and only if
$|L'/s^{-1}\Z^2|=p$. Therefore from (\ref{esymm1}) we get
\begin{equation} \label{eheck1}
\pi(v_p)\delta_L=\sum_{L'\supset L\colon |L'/L|=p}\delta_{L'}.
\end{equation}
We also have
\begin{equation} \label{eheck2}
\pi(u_p)\delta_L=\delta_{p^{-1}L}.
\end{equation}
These are classical formulas for Hecke operators.

Note that the algebras $\hecke{S_p}{\Gamma}$ and
$C(\G(\Z_p)\backslash\mtwo(\Z_p))\subset
C(\G(\rr)\backslash\mtwo(\rr))$ are subalgebras of
$\mathfrak{A}_p^{\Z^*_p}={\redheck{\mtwo(\Z[p^{-1}]) \rtimes
G_p}{\mz\rtimes\Gamma}}^{\Z^*_p}$. Recall that equation
\eqref{esymm2} defines a representation $\pi_p$ of $\mathfrak{A}_p$
on $\ell^2(\Gamma\backslash S_p)$. We can now establish the
following.

\begin{proposition} \label{Irr}
For each prime $p$, the C$^*$-algebra generated by
$\hecke{S_p}{\Gamma}$ and $C(\G(\Z_p)\backslash\mtwo(\Z_p))$ in the
representation $\pi_p$ contains the algebra of compact operators on
$\ell^2(\Gamma\backslash S_p)$. In particular, the restriction
of~$\pi_p$ to $\mathfrak{A}_p^{\Z^*_p}$ is irreducible.
\end{proposition}

\bp First of all note that for the representation $\pi_p$ both
formulas  (\ref{eheck1})-(\ref{eheck2}) are valid. Given
$L=s^{-1}\Z^2$ with $s\in S_p$, we denote by $e_L\in
C(\G(\Z_p)\backslash\mtwo(\Z_p))$ the characteristic function of the
set $\G(\Z_p)s$, and view it as an element of $C(\G(\rr) \backslash
\mtwo(\rr))$. Then
$$
e_{L}=\pi_L-\bigvee_{L'\supset L\colon
|L'/L|=p}\pi_{L'}=\prod_{L'\supset L\colon |L'/L|=p}
(\pi_L-\pi_{L'}).
$$
Since $\pi_p(e_L)(\delta_{L^{''}})$ can be nonzero only if $L\subset
L^{''}\subsetneqq L'$, and $\vert L'/L\vert=p$, it follows that
$\pi_p(e_L)$ is the projection onto
$\C\delta_L\subset\ell^2(\Gamma\backslash S_p)$.

Now suppose that $\det(s)=p^n$. Then by (\ref{eheck1}) the vector
$\pi_p(v^n_p)\delta_{\Z^2}$ is a linear combination with nonzero
coefficients of the vectors $\delta_{L'}$, where $L'$ runs over all
lattices such that $[L'\colon\Z^2]=p^n$. It follows that
$\pi_p(e_Lv^n_p)\delta_{\Z^2}=\lambda\delta_L$ for a nonzero scalar
$\lambda$. Then $\lambda^{-1}\pi_p(e_Lv^n_pe_{\Z^2})$ is a partial
isometry with initial space $\C\delta_{\Z^2}$ and range
$\C\delta_L$, and the proposition follows. \ep

It seems natural to refer to the  C$^*$-subalgebra of
$\mathfrak{A}_p^{\Z^*_p}$ generated by $\hecke{S_p}{\Gamma}$ as
\emph{the Toeplitz-Hecke algebra at prime $p$}. The representation
$\pi_p$ restricted to this algebra is no longer irreducible for
instance because the operators in the image commute with the action
of $\Gamma$ on $\ell^2(\Gamma\backslash S_p)$ defined by
$\gamma\delta_L=\delta_{\gamma L}$. In the rest of the section we
shall show that nevertheless the image under $\pi_p$ of the
Toeplitz-Hecke algebra at prime $p$ contains the projection onto
$\C\delta_{\Z^2}$. We begin with a lemma.

\begin{lemma}\label{need_a_lemma}
Assume $L=s^{-1}\Z^2$ for $s\in S_p$ is such that $L\ne\Z^2$ and $L$
does not contain $p^{-1}\Z^2$. Then \enu{i} there is a unique
lattice $L'\supset L$ such that $[L'\colon L]=p$ and
$p^{-1}\Z^2\subset L'$; \enu{ii} there is a unique lattice
$L''\subset L$ such that $[L\colon L'']=p$ and $\Z^2\subset
L''$.
\end{lemma}

\bp Since $L/\Z^2$ is a nontrivial $p$-group, it has elements of
order $p$. In other words, $L\cap p^{-1}\Z^2$ strictly contains
$\Z^2$ and is strictly contained in $p^{-1}\Z^2$. Hence $(L\cap
p^{-1}\Z^2)/\Z^2$ is a group of order $p$. Then $L'=L+p^{-1}\Z^2$
has the properties $[L'\colon L]= [p^{-1}\Z^2\colon L\cap
p^{-1}\Z^2]=p$ and $L'\supset p^{-1}\Z^2$.  If there exists another
lattice $L''$ with such properties then $L=L'\cap L''$ and hence
$p^{-1}\Z^2\subset L$, which is a contradiction. Thus (i) is proved.

\smallskip

Since $L/\Z^2$ is a nontrivial $p$-group, it has a subgroup of index
$p$. The preimage of such a subgroup in $L$ is a lattice with the
properties required in (ii). Moreover, we have a one-to-one
correspondence between such lattices and subgroups of $L/\Z^2$ of
index $p$. Since any double coset in $S_p$ has a diagonal
representative, $L/\Z^2$ has the form $\Z/p^a\Z\oplus\Z/p^b\Z$,
$a\le b$. Then the condition that $L$ does not contain~$p^{-1}\Z^2$
means exactly that $a=0$, or equivalently, $L/\Z^2$ is a cyclic
$p$-group (indeed, to say that $L$ does not contain $p^{-1}\Z^2$ is
the same as saying that the group of elements $x\in L/\Z^2$ such
that $px=0$ contains strictly fewer than $p^2$  elements). But if
$L/\Z^2$ is cyclic then it contains only one subgroup of index $p$.
Thus (ii) is also proved. \ep

We can now prove our claim about the projection.

\begin{proposition} \label{projection}
The image of $v^*_pv_p-v_pv^*_p-p(1-u_pu^*_p)$ under $\pi_p$ is the
orthogonal projection onto~
$\C\delta_{\Z^2}$  in $B(\ell^2(\Gamma\backslash S_p))$.
\end{proposition}

\bp Let $T:=v^*_pv_p-v_pv^*_p-p(1-u_pu^*_p)$. By virtue of
(\ref{eheck1}) we have
$$
\pi_p(v^*_p)\delta_L=\sum_{L''\subset L\colon\Z^2\subset
L'',\,|L/L''| =p}\delta_{L''}.
$$
Since any lattice $L''$ containing a lattice $L$ as a subgroup of
index $p$ is contained in $p^{-1}L$, and any sublattice of $p^{-1}L$
of index $p$ contains $L$, we get by faithfulness of $\pi_p$ that
$v_p^*u_p=v_p$. We can then write
$$
v^*_pv_p-v_pv^*_p=v^*_p(1-u_pu_p^*)v_p.
$$

To compute the action of $\pi_p(T)$ on
$\delta_L\in\ell^2(\Gamma\backslash S_p)$, assume first that $L$
contains $p^{-1}\Z^2$. Since the operator $1-\pi_p(u_pu^*_p)$ is the
projection onto the space spanned by the $\delta_{L'}$'s such that
$L'$ does not contain $p^{-1}\Z^2$,
$$
\pi_p((1-u_pu_p^*)v_p)\delta_L=0=\pi_p(1-u_pu_p^*)\delta_L.
$$
Thus $\pi_p(T)\delta_L=0$. Assume now that $L$ does not contain
$p^{-1}\Z^2$. Since $R_\Gamma\diag{1}{p}=p+1$, see e.g.
\cite[Chapter~IV]{kri}, there exist exactly $p+1$ lattices
containing $L$ as a subgroup of index $p$, and by
Lemma~\ref{need_a_lemma}(i) only one of them contains $p^{-1}\Z^2$.
Let $L_1,\dots,L_p$ denote the remaining $p$ lattices. Then
$$
\pi_p((1-u_pu_p^*)v_p)\delta_L=\sum^p_{i=1}\delta_{L_i}.
$$
By Lemma~\ref{need_a_lemma}(ii) applied to $L_i$ for each $i$, the
lattice $L$ is the unique lattice with the properties that
$[L_i\colon L]=p$ and $\Z^2\subset L$. Thus
$\pi_p(v^*_p)\delta_{L_i}=\delta_L$, and hence
$$
\pi_p(v^*_p(1-u_pu_p^*)v_p)\delta_L
=p\delta_L=p\pi_p(1-u_pu_p^*)\delta_L,
$$
giving $\pi_p(T)\delta_L=0$. For $L=\Z^2$ the computation is
similar, but now all $p+1$ lattices containing $\Z^2$ as a
sublattice of index~$p$ do not contain $p^{-1}\Z^2$. Hence
$$
\pi_p(v^*_p(1-u_pu_p^*)v_p)\delta_{\Z^2}=(p+1)\delta_{\Z^2},
$$
and the result follows. \ep

\bigskip

\section{KMS-states} \label{SectionKMS}

In view of \proref{MainHecke}, the canonical dynamics, as defined in
\cite[Proposition 4]{bos-con} for a general Hecke algebra, is
determined by  $\sigma_t(f)(g,m) =\det(g)^{2it}f(g,m)$ on the finite
part of the Connes-Marcolli system, but for simplicity we shall omit
the factor~$2$ and consider instead the dynamics
$$
\sigma_t(f)(g,m)=\det(g)^{it}f(g,m)
$$
for every $f$ in the dense subalgebra $
C_c(\fibb{\glap}{\slr}{\mr})$ of $\redheck{P}{P_0}$.

\subsection{KMS-states on corners of crossed products}
Before we turn our attention to the classification of KMS-states on
$\redheck{P}{P_0}$, we first prove some general results on
KMS-states on crossed products, which will be needed later, and may
prove useful elsewhere.

Suppose that $X$ is a locally compact space, $X_0$ a clopen subset
of $X$, $G$ a locally compact group acting on $X$, and $G_0$ a
compact open subgroup of $G$. Assume $G_0X_0=X_0$ and $\bigcup_{g\in
G}gX_0=X$. Let $p_{G_0}=\int_{G_0}\lambda_g\,dg$, and denote by $p$
the projection $\31_{X_0}p_{G_0}$ in the multiplier algebra of
$C_0(X)\rtimes_r G$.

For the purpose of this subsection we have to assume that the action of $G$ or $G_0$ is free. Although this assumption is not satisfied in our main example, it will be satisfied for certain subsystems, see the proof of Lemma~\ref{KMS1}.

\begin{lemma}\label{lemma_p_full}
Suppose that the action of $G_0$ on $X_0$ is free. Then the
projection $p$ is full.
\end{lemma}

\bp It is well-known that the projection $p_{G_0}$ is full in the
multiplier algebra of $C_0(X_0)\rtimes_r G_0$ (this is essentially
equivalent to the fact that $C_0(X_0)\rtimes_r G_0$ and
$C_0(G_0\backslash X_0)$ are Morita equivalent), but for
completeness we provide a short argument. The ideal generated by
$p_{G_0}$ is the closed linear span of functions of the form
$(g,x)\mapsto f_1(x)f_2(gx)$, $f_i\in C_0(X_0)$. Since by the
assumption on freeness such functions separate points, and thus span
a dense subspace of $C_0(G_0\times X_0)$, it follows that $p_{G_0}$
is full in $C_0(X_0)\rtimes_r G_0$.

Hence the ideal generated by $p$ in $C_0(X)\rtimes_r G$ contains
$C_0(X_0)\rtimes_r G_0$. Since $\bigcup_{g\in G}gX_0=X$, the
elements $\lambda_{g_1}a\lambda_{g_2}$ for  $a\in C_0(X_0)\rtimes_r
G_0$ and $g_i\in G$ span a dense subspace of $C_0(X)\rtimes_r G$.
Thus $p$ is full. \ep

Assume now that $N\colon G\to \R^*_+$ is a homomorphism whose kernel
contains $G_0$. Define a dynamics~$\sigma$ on $C_0(X)\rtimes_r G$ by
$\sigma_t(f)(g,x)=N(g)^{it}f(g,x)$. We want to classify
$\sigma$-KMS$_\beta$-states on
$$
A:=p(C_0(X)\rtimes_r G)p.
$$

Recall, see e.g.~\cite{kus0}, that a semifinite $\sigma$-invariant
weight $\varphi$ is called a $\sigma$-KMS$_\beta$-weight, where
$\beta\in\R$, if
$$
\varphi(aa^*)=\varphi(\sigma_{i\beta/2}(a)^*\sigma_{i\beta/2}(a))
$$
for all $\sigma$-analytic elements $a\in A$. Since the projection
$p$ is full, any KMS$_\beta$-state on $A$ extends uniquely to a
KMS$_\beta$-weight on $C_0(X)\rtimes_r G$, see
e.g.~\cite[Remark~3.3(i)]{LN}.

On the other hand, if we let $E$ be the C$^*$-valued weight
$$
E\colon C_0(X)\rtimes_r G\to C_0(X)
$$
such that $E(f)=f(e,\cdot)$, then any Radon measure $\mu$ on~$X$
defines a semifinite weight $\mu_*\circ E$ on $C_0(X)\rtimes_r G$.
This is a standard way of getting dual weights on von Neumann
algebras~\cite{str}, and can be justified in the C$^*$-algebra
setting using~\cite{kus}. We are not going into details because
there are at least two other ways of constructing this weight:
either by restricting the dual weight on the von Neumann algebra
crossed product to the C$^*$-algebra crossed product, or by
using~\cite{qu-ve}.

Note that the equation $E_0(a)=E(a)p$ defines a conditional
expectation $A\to C_0(G_0\backslash X_0)p\subset A$, and by
identifying $C_0(G_0\backslash X_0)p$ with $C(G_0\backslash X_0)$ we
get
\begin{equation}\label{small_condexp}
(\mu_*\circ E)\vert_A=(\mu|_{X_0})_*\circ E_0.
\end{equation}

\begin{theorem} \label{KMSf}
Assume that the action of $G$ on $X$ is free, $G_0X_0=X_0$ and
$\bigcup_{g\in G}gX_0=X$. Then for each $\beta \in \R$ there is a
one-to-one correspondence between $\sigma$-KMS$_\beta$-states on
$A=p(C_0(X)\rtimes_r G)p$ and Radon measures $\mu$ on $X$ such that
\begin{equation}\label{scaling_Radon_meas}
\mu(X_0)=1\ \ \text{and}\ \ \mu(gY)=N(g)^{-\beta}\mu(Y)\text{ for
compact } Y\subset X \text{  and  }g\in G.
\end{equation}
\end{theorem}

\begin{proof}
A standard computation  using the covariance relation in $A$  shows
that any measure satisfying the conditions in
\eqref{scaling_Radon_meas} defines a KMS$_\beta$-state. This is in
fact valid without the assumption that $G$ acts freely on $X$. We
omit the details.

The nontrivial part is to show that any KMS-state is determined by a
measure. Assume $\varphi$ is a $\sigma$-KMS$_\beta$-state on~$A$.
Using \cite[Remark~3.3(i)]{LN} we extend it to a
$\sigma$-KMS$_\beta$-weight on $C_0(X)\rtimes_r G$, which we
continue to denote by $\varphi$. Consider the subsets
$$
\mathfrak N_\varphi=\{a\,|\,\varphi(a^*a)<\infty\},\ \ \mathfrak
M_\varphi=\lspa\,\mathfrak N_\varphi^*\mathfrak
N_\varphi=\lspa\{a\ge0\,|\,\varphi(a)<\infty\}
$$
of $C_0(X)\rtimes_r G$. To continue, recall from e.g.~\cite{kus0},
\cite{str} that $\varphi$ extends uniquely to a linear functional
on~$\mathfrak M_\varphi$, and $\mathfrak M_\varphi$ is a bimodule
over the algebra of $\sigma$-analytic elements. Since $A\subset
\mathfrak M_\varphi$, and $f\in C_0(X)$ and $\lambda_g$ are
analytic, we conclude as in the proof of Lemma~\ref{lemma_p_full}
that $\mathfrak M_\varphi\cap (C_0(X)\rtimes G_0)$ is dense in
$C_0(X)\rtimes G_0$. Hence the restriction of $\varphi$ to
$C_0(X)\rtimes G_0$ is a semifinite weight. Because $G_0\subset \ker
N$, this weight is in fact a trace on $C_0(X)\rtimes G_0$. Since
$p_{G_0}$ is full and $p_{G_0}(C_0(X)\rtimes
G_0)p_{G_0}=C_0(G_0\backslash X)p_{G_0}$, this trace is uniquely
determined by a weight on $C_0(G_0\backslash X)$, that is, by a
$G_0$-invariant measure $\mu$ on $X$. But $(\mu_*\circ
E)|_{C_0(X)\rtimes G_0}$ is also a semifinite trace on
$C_0(X)\rtimes G_0$ which extends the weight
$fp_{G_0}\mapsto\mu_*(f)$ on $C_0(G_0\backslash X)p_{G_0}$, and so
we conclude that
$$
\varphi=\mu_*\circ E\ \ \text{on}\ \ C_0(X)\rtimes G_0.
$$
In particular, $f\in\mathfrak M_\varphi$ for $f\in C_c(G_0\times
X)$, and $\varphi(f)=\int_Xf(e,\cdot)d\mu$. Using again that
$\lambda_g$ is $\sigma$-analytic for every $g\in G$, we conclude
that $f\in\mathfrak M_\varphi$ for $f\in C_c(G\times X)$. If the
support of $f$ does not intersect $\{e\}\times X$ then, since the
action is free, we can find functions $f_1,\dots,f_n\in C_c(X)$ such
that if $(g,x)$ is in the support of $f$ for some $g\in G\setminus
\{e\}$ and $x\in X$, then $\sum^n_{i=1}f_i(x)=1$ and $f_i(gx)=0$ for
all $i$. It follows that $f_i*f=0$, so by using the KMS-condition we
get
$$
\varphi(f*f_i)=\varphi(f_i*f)=0\ \text{for}\ i=1, \dots ,n.
$$
Hence $\varphi(f)=0$, and because $G_0$ is open in $G$, we conclude
that
\begin{equation}\label{phi_on_big_algebra}
\varphi(f)=\int_Xf(e,\cdot)d\mu
\end{equation}
for $f\in C_c(G\times X)$. Now it is easy to check that $\mu$
satisfies the conditions in \eqref{scaling_Radon_meas}, for example
$\mu(X_0)=1$ because $\varphi(p)=1$.

The equality \eqref{phi_on_big_algebra} completely determines the
state $\varphi$ on $A$. Since $\mu$ satisfies the scaling condition
in \eqref{scaling_Radon_meas}, we can also conclude that $\mu_*\circ
E$ is a $\sigma$-KMS$_\beta$-weight on $C_0(X)\rtimes_r G$ which
extends $\varphi$ on~$A$. Since the extension is unique, we must
have $\varphi=\mu_*\circ E$ on $C_0(X)\rtimes_r G$. \ep

\subsection{Classification of KMS-states on $\redheck{P}{P_0}$}
We can not apply Theorem~\ref{KMSf} directly to $X=\ma$, $X_0=\mr$,
$G=\glap$, $G_0=\slr$ and $N(g)=\det(g)$, as we would like to,
because the action is very far from being free. So we shall impose
an additional assumption on KMS-states, and apply the theorem to
systems corresponding to finite sets of primes.

Let us first recall that for every $\beta>1$ there exists a
canonical measure $\mu_{\beta,f}$ on $\ma$
satisfying~\eqref{scaling_Radon_meas}. As we already remarked in the
proof of Theorem~\ref{KMSf}, any such measure defines a
KMS$_\beta$-state, but in view of lack of freeness we can not be
sure that in this way we get all KMS-states.

The construction of $\mu_{\beta,f}$ is as follows,
see~\cite[Section~4]{LLN} for details. For each prime number $p$
consider the Haar measure on $\G(\Z_p)$ normalized so  that the
total mass is $(1-p^{-\beta})(1-p^{-\beta+1})$. This measure extends
to a unique measure $\mu_{\beta,p}$ on $\G(\Q_p)$ such that
$$
\mu_{\beta,p}(gZ)=|\det(g)|_p^\beta\mu_{\beta,p}(Z)\ \ \text{for
compact}\ \ Z\subset\G(\Q_p)\ \ \text{and}\ \ g\in\G(\Q_p),
$$
where $\vert a\vert_p$ denotes the $p$-adic valuation of $a$. The
total mass of $\mtwo^i(\Z_p):=\mtwo(\Z_p)\cap\G(\Q_p)$, which is the
set of regular matrices, is one, and  we can define a measure on
$\mtwo(\af)$ by setting $\mu_{\beta,f}=\prod_p\mu_{\beta,p}$. By
construction this measure satisfies $\mu_{\beta,f}(\mtwo(\rr))=1$
and
$$
\mu_{\beta,f}(gZr)=\left(\prod_p|\det(g_p)\vert_p\right)^\beta
\mu_{\beta,f}(Z)
$$
for $Z\subset\mtwo(\af)$, $g\in\G(\af)$ and $r\in\G(\rr)$. Note that
it is not difficult to show that $\mu_{\beta,f}(Zg)=(\prod_p
\vert\det(g_p)\vert_p)^\beta\mu_{\beta,f}(Z)$ for $g\in\gla$, but we
will not need this. Note also that $\mu_{2,f}$ is a Haar measure on
$\mtwo(\af)$. Furthermore, by construction of $\mu_{\beta,f}$ the
set
\begin{equation}\label{def_non_singular}
\mtwo^i(\rr):=\prod_{p\in\primes}\mtwo^i(\Z_p)
\end{equation}
is a subset of $\mtwo(\rr)$ of full measure. We denote by
$\varphi_\beta$ the  KMS$_\beta$-state corresponding to~$\mu_{\beta,
f}$ for $\beta >1$.

On the other hand, for $\beta>2$ and every $w\in\G(\rr)$ we can
construct a $\sigma$-KMS$_\beta$-state as follows. Consider the
representation $\pi_w$ of $\redheck{P}{P_0}$ on
$\ell^2(\Gamma\backslash S)$ introduced in
Theorem~\ref{pi_faithful_our_pair}. Define an unbounded positive
selfadjoint operator $H$ on $\ell^2(\Gamma\backslash S)$ by
\begin{equation} \label{egen}
H\delta_L=\log[L\colon\Z^2]\delta_L.
\end{equation}
Then, see e.g. \cite[Lemma~1.18]{con-mar}, we have
\begin{equation} \label{ezeta}
\Tr(e^{-\beta H})=\zeta(\beta)\zeta(\beta-1),
\end{equation}
where $\zeta$ is the Riemann $\zeta$-function. The dynamics $\sigma$
is implemented in the representation $\pi_w$ by the one-parameter
unitary group with generator $H$, that is,
$$
\pi_w(\sigma_t(a))=e^{itH}\pi_w(a)e^{-itH}.
$$

\begin{lemma}\label{eextr} $($cf. \cite[Section 1.7]{con-mar}$)$
The formula
$$
\varphi_{\beta,w}(a)=\zeta(\beta)^{-1}\zeta(\beta-1)^{-1}
\Tr(\pi_w(a)e^{-\beta H})
$$
defines a $\sigma$-KMS$_\beta$-state that depends only on
$\det(w)\in\rr^*$.
\end{lemma}

\bp Given $w,w'\in \G(\rr)$ such that $\det(w)=\det(w')$, similarly
to Remark~\ref{pi_w_equiv} we conclude that the unitary operator~$U$
which maps $\delta_L$ onto $\delta_{w'w^{-1}L}$ has the property
$U\pi_w(\cdot)U^*=\pi_{w'}$. Since $U$ commutes with $H$, we
conclude that $\varphi_{\beta,w}=\varphi_{\beta,w'}$. \ep

We shall classify the following class of KMS-states: given a
KMS-state $\varphi$ on $\redheck{P}{P_0}$, we restrict it to the
subalgebra $C(\slr\backslash\mtwo(\rr))$,  obtaining  an
$\slr$-invariant probability measure~$\mu_\varphi$ on $\mtwo(\rr)$.
We say that $\varphi$ is \emph{regular at every prime} if the set
$\mtwo^i(\rr)$ from \eqref{def_non_singular} is a subset of full
measure for $\mu_{\varphi}$. Equivalently, for each $p$ in
$\primes$, the push-forward of $\mu_\varphi$ under the projection
$\mtwo(\rr)\to\mtwo(\Z_p)$ gives a measure on $\mtwo(\Z_p)$ for
which the set of singular matrices $\mtwo(\Z_p)\setminus
\mtwo^i(\Z_p)$ has measure zero. Apart from the fact that then
$\varphi$ is in some sense supported on regular matrices, another
reason to call $\varphi$ regular is the following result.

\begin{lemma}\label{normal_state}
A $\sigma$-KMS$_\beta$-state $\varphi$ on $\redheck{P}{P_0}$ is
regular at every prime if and only if its restriction to
$\mathfrak{A}_p^{\Z^*_p}={\redheck{\mtwo(\Z[p^{-1}]) \rtimes
G_p}{\mz\rtimes\Gamma}}^{\Z^*_p}$ is normal with respect to the
representation $\pi_p$ defined by \eqref{esymm2} for every~$p$ (in
other words, since the representation is irreducible by
Proposition~\ref{Irr}, $\varphi$ extends to a normal state on
$B(\ell^2(\Gamma\backslash S_p))$).
\end{lemma}

\bp We have $\mtwo^i(\Z_p)=\bigcup_{s\in S_p} \G(\Z_p)s$, which can
be established similarly to \eqref{decomp_of_glap}. Recall from the
proof of Proposition~\ref{Irr} that for a prime $p$ and a lattice
$L=s^{-1}\Z^2$ with $s\in S_p$, we denoted by $e_L$ the
characteristic function of $\G(\Z_p)s$. Then to say that $\varphi$
is regular is the same as requiring that
$$
\sum_{L\supset\Z^2\colon L/\Z^2\ \text{is a}\ p\text{-group}}
\varphi(e_L)=1.
$$
But since $\pi_p(e_L)$ is the projection onto $\C\delta_L$, this
condition just means that the restriction of $\varphi$ to the
algebra of compact operators on $\ell^2(\Gamma\backslash S_p)$ is a
state. It remains to recall that if we have a state on a
C$^*$-algebra $A\subset B(H)$ containing the algebra $K(H)$ of
compact operators then this state extends to a normal state on
$B(H)$ if and only if its restriction to $K(H)$ is a state. \ep

We next prove that a regular KMS-state is completely determined
by $\mu_\varphi$ on $\mtwo(\rr)$.

\begin{lemma} \label{KMS1}
For all $\beta\in \R$ there is a one-to-one correspondence between
$\sigma$-KMS$_\beta$-states on $\redheck{P}{P_0}$ which are regular
at every prime, and measures $\mu$ on $\ma$ such that
$\mu(\mtwo(\rr))=\mu(\mtwo^i(\rr))=1$ and
$\mu(gY)=\det(g)^{-\beta}\mu(Y)$ for compact $Y\subset \ma$ and
$g\in\glap$.
\end{lemma}

\bp Certainly a measure $\mu$ with the properties described in the
lemma gives rise to a KMS-state which is regular.

Conversely, suppose that $\varphi$ is a regular KMS-state. We then
have a measure $\mu_\varphi$ on $\mtwo(\rr)$ such that
$\mu_\varphi(\mtwo^i(\rr))=1$, and to get the extra properties of
$\mu_\varphi$ we shall use the push-forwards of $\mu_\varphi$ under
the projection maps from $\mtwo(\rr)$ onto the coordinates
corresponding to finite subsets of primes, whereby we will be in a
position to apply Theorem~\ref{KMSf}.

For each finite subset $F$ of primes denote by $\Q_F$ the ring
$\prod_{p\in F}\Q_p$, and by $\Z_F$ its subring $\prod_{p\in
F}\Z_p$. Consider the subgroup $G_F$ of $\glq$ generated by the
subgroups $G_p$ for $p\in F$. Denote by $\glp(\Q_F)$ the closure of
$G_F$ in $\G(\Q_F)$. Similarly to Lemma~\ref{glqdenseinglap} one can
show that $\glp(\Q_F)$ equals $G_F \sltwo(\Z_F)$, and is the group
of elements $g\in\G(\Q_F)$ such that $\det(g)$ is  an element of the
multiplicative group generated by elements $p\in F$ considered as a
subgroup of $\Q_F$. Then $\glp(\Q_F)$ acts on $\mtwo(\Q_F)$ by
multiplication on the left, and we define
$$
\mathfrak{A}_F=p_F(C_0(\mtwo(\Q_F))\rtimes_r\glp(\Q_F))p_F,
$$
where $p_F$ is the projection in the crossed product corresponding
to the characteristic function of the compact open subset
$\sltwo(\Z_F)\times\mtwo(\Z_F)$ of $\glp(\Q_F)\times \mtwo(\Q_F)$.
We view $\mathfrak{A}_F$ as a subalgebra of
$p_0(C_0(\ma)\rtimes_r\glap)p_0$ by identifying it with the closure
of ($\slr\times\slr$)-invariant functions $f$ on
$$
\{(g,m)\,\mid\,g\in\glap,\ m\in\mtwo(\rr),\ gm\in\mtwo(\rr)\}
$$
such that $f(g,m)=0$ if $g\notin G_F\slr$, and $f(g,m)=f(g,m')$ if
$m_p=m'_p$ for all $p\in F$. Alternatively, with $N_F:=\prod_{p\in
F}p$, one can show that $\mathfrak{A}_F$ is the reduced Hecke
C$^*$-algebra of the pair
$$
\left(\matr{\mtwo(\Z[N_F^{-1}])}{\glp(\Z[N_F^{-1}])},
\matr{\mz}{\slz} \right).
$$

The restriction of $\varphi$ to $\mathfrak{A}_F$ is a KMS-state,
which upon further restriction to the subalgebra
$C(\sltwo(\Z_F)\backslash\mtwo(\Z_F))$ defines an
$\sltwo(\Z_F)$-invariant measure $\mu_F$ on $\mtwo(\Z_F)$. Note that
$\mu_F$ is the push-forward of $\mu_\varphi$ to $\mtwo(\Z_F)$. Thus,
since $\varphi$ is regular, $\mu_F(\mtwo(\Z_F)\setminus \prod_{p\in
F}\mtwo^i(\Z_p))=0$, and therefore the restriction of $\varphi$ to
the ideal
$$
{I}_F=p_F(C_0(\G(\Q_F))\rtimes_r\glp(\Q_F))p_F
$$
is still a state. In particular, $\varphi|_{\mathfrak{A}_F}$ is
completely determined by $\varphi|_{{I}_F}$. Since $\glp(\Q_F)$
certainly acts freely on $\G(\Q_F)$,  Theorem~\ref{KMSf} applied to
the algebra ${I}_F$ says that the measure $\mu_F$ uniquely extends
to a measure on $\mtwo(\Q_F)$, which we still denote by $\mu_F$,
such that $\G(\Q_F)$ is a subset of $\mtwo(\Q_F)$ of full measure,
$\mu_F(\mtwo(\Z_F))=1$, and
$$
\mu_F(gY)=\det(g)^{-\beta}\mu_F(Y)
$$
for $g\in\glp(\Q_F)$ and $Y\subset\mtwo(\Q_F)$ compact.

If $F'$ is another finite set of prime numbers containing $F$, then
the push-forward of $\mu_{F'}$ under the projection map
$$
\mtwo(\Q_{F'})\supset\mtwo(\Q_F)\times\mtwo(\Z_{F'\setminus F})
\to\mtwo(\Q_F)
$$
must coincide with $\mu_F$. It follows that there exists a unique
measure $\mu$ on $\ma$ such that its push-forward under
$$
\ma\supset\mtwo(\Q_F)\times\prod_{q\notin F}\mtwo(\Z_q)
\to\mtwo(\Q_F)
$$
coincides with $\mu_F$.  By the properties of the $\mu_F$'s we
conclude that $\mu$ is $\slr$-invariant, the scaling condition
$\mu(gY)=\det(g)^{-\beta}\mu(Y)$ holds for compact $Y\subset \ma$
and $g\in\glq$, and that~$\mu(\mtwo(\rr))=\mu(\mtwo^i(\rr)) =1$.

Since $\bigcup_F\mathfrak{A}_F$ is dense in $\redheck{P}{P_0}$, the
state $\varphi$ is completely determined by $\mu$: indeed, $\varphi$
is obtained by composing $\mu|_{\mtwo(\rr)}$ with the conditional
expectation $E_0:\redheck{P}{P_0}\to C(\slr\backslash\mtwo(\rr))$
from \eqref{small_condexp}, because this is true on each
$\mathfrak{A}_F$, i.e. $\varphi\vert_{\mathfrak{A}_F}$ is the
composition of $\mu_F$ with the similar conditional expectation onto
$C(\sltwo(\Z_F)\backslash \mtwo(\Z_F))$. Therefore the proof is
complete. \ep

It is not difficult to show that for $\beta\ne0,1$ the condition
$\mu(\mtwo^i(\rr))=1$ follows from the remaining ones. In fact, a
stronger result is proved in~\cite[Corollary~3.6]{LLN}, where it is
shown that a weaker condition of relative invariance under Hecke
operators implies that the measure of singular matrices must be
zero.  On the other hand, for $\beta=0,1$ one does have measures
satisfying all conditions of Lemma~\ref{KMS1} except
$\mu(\mtwo^i(\rr))=1$, see~\cite[Remark~4.7]{LLN}.

\smallskip

We can now formulate our main classification result.

\begin{theorem}\label{main_theorem}
For $\beta\in\R$ denote by $K_\beta$ the simplex of
$\sigma$-KMS$_\beta$-states, and by $K'_\beta$ the subset of states
which are regular at every prime. Then $K'_\beta$ is a subsimplex of
$K_\beta$, that is, it is closed and if a probability measure on
$K_\beta$ has barycenter in $K'_\beta$ then the measure must be
supported on $K'_\beta$. Furthermore, we have: \enu{i} for
$\beta\le1$ the set $K'_\beta$ is empty; \enu{ii} for
$\beta\in(1,2]$ the set $K'_\beta$  consists of one point
$\varphi_\beta$ corresponding to the measure $\mu_{\beta,f}$;
\enu{iii}\, for $\beta>2$ the simplex $K'_\beta$ is isomorphic to
the simplex of probability measures on $\slr\backslash\glr$; in
particular, extremal points of $K'_\beta$ correspond to
$\slr$-orbits in $\glr$, and the state corresponding to $\slr w$ is
$\varphi_{\beta,w}$.
\end{theorem}

\bp For every prime $p$ consider the operator $H_p$ on
$\ell^2(\Gamma\backslash S_p)$ defined exactly as $H$ in
(\ref{egen}). Then similarly to \eqref{ezeta} we have
\begin{equation} \label{ezetap}
\Tr(e^{-\beta H_p})=\begin{cases}
(1-p^{-\beta})^{-1}(1-p^{-\beta+1})^{-1},&
\text{if}\ \beta>1,\\
+\infty,&\text{otherwise}.
\end{cases}
\end{equation}

The operator $H_p$ implements the dynamics $\sigma$ on
$\mathfrak{A}_p$ in the representation $\pi_p$, in the sense that
$\pi_p(\sigma_t(a))=e^{itH_p}\pi_p(a)e^{-itH_p}$. It follows that if
a $\sigma$-KMS$_\beta$-state on $\mathfrak{A}_p^{\Z_p^*}$ extends to
a normal state on $B(\ell^2(\Gamma\backslash S_p))$, then this
extension is a KMS$_\beta$-state for the dynamics $\Ad e^{itH_p}$.
By virtue of (\ref{ezetap}) the latter dynamics admits no normal
KMS-states for $\beta\le1$, while for $\beta>1$ there exists a
unique normal KMS-state on $B(\ell^2(\Gamma\backslash S_p))$ given
by
$$
\psi_{\beta,p}=(1-p^{-\beta})(1-p^{-\beta+1}) \Tr(\cdot\,
e^{-\beta H_p}).
$$

Thus Lemma~\ref{normal_state} implies that for $\beta>1$, a
$\sigma$-KMS$_\beta$-state $\varphi$ is regular at every prime if
and only if $\varphi=\psi_{\beta,p}\circ\pi_p$ on
$\mathfrak{A}_p^{\Z_p^*}$ for every $p$, while for $\beta\le1$ there
are no regular KMS-states. This  proves (i) and shows that
$K'_\beta$ is closed. Furthermore, assuming $\beta>1$, every
KMS-state on $\redheck{P}{P_0}$ defines a positive KMS-functional on
the algebra $K(\ell^2(\Gamma\backslash S_p))$ of compact operators.
Since the latter algebra has a unique KMS-state, this functional is
$\lambda\psi_{\beta,p}$ for some $\lambda\le1$. Moreover, since
$\psi_{\beta,p}$ is faithful, to check whether $\lambda=1$ it
suffices to evaluate the functional on one positive nonzero
operator. In other words, if we fix a positive nonzero element
$a_p\in\mathfrak{A}_p^{\Z^*_p}$ such that $\pi_p(a_p)$ is compact,
then $\varphi(a_p)\le\psi_{\beta,p}(\pi_p(a_p))$, and the equality
holds for all $p$ if and only if $\varphi$ is regular at every
prime. It is now clear that if the barycenter of a probability
measure on the set of KMS-states is a regular state then almost
every state is regular.

\smallskip

It remains to prove (ii) and (iii). In other words, we want to
classify all measures $\mu$ satisfying the conditions in
Lemma~\ref{KMS1}. Since similar results are proved in \cite{LLN}, we
will here be somewhat brief.

\smallskip

Assume $1<\beta\le2$. Then by the proof of \cite[Lemma~4.5]{LLN},
for every measure $\mu$ on $\ma$ such that $\mu(\mtwo(\rr))=1$ and
$\mu(gY)=\det(g)^{-\beta}\mu(Y)$ for $Y\subset\ma$ and $g\in\glq$,
the action of $\glq$ on $(\ma,\mu)$ is ergodic. Since the set of
such measures is convex, this means that this set consists of at
most one point. Hence $\mu_{\beta,f}$ is the unique such measure.

\smallskip

Assume now that $\beta>2$, and let $\mu$ be a measure on $\ma$
satisfying the conditions of Lemma~\ref{KMS1}. For a finite set $F$
of prime numbers denote
$$
Y_F:=\G(\Z_F)\times\prod_{q\notin F}\mtwo(\Z_q),
$$
and $S_F:=G_F\cap\mzp$. Then by regularity the set $S_FY_F$ is a
subset of $\mtwo(\rr)$ of full measure. Therefore
$$
1=\mu(S_FY_F)=\sum_{s\in S_F/\Gamma}\mu(sY_F)=\mu(Y_F) \sum_{s\in
S_F/\Gamma}\det(s)^{-\beta}
=\mu(Y_F)\prod_{p\in
F}(1-p^{-\beta})^{-1}(1-p^{-\beta+1})^{-1},
$$
where the last step follows e.g. from \cite[Equation (3.3)]{LLN}.
Thus $\mu(Y_F)=\prod_{p\in F}(1-p^{-\beta})(1-p^{-\beta+1})$. By
taking the intersection of all sets $Y_F$ for different $F$ we
obtain
$$
\mu(\G(\rr))=\prod_{p\in\primes}(1-p^{-\beta})(1-p^{-\beta+1})
=\zeta(\beta)^{-1}\zeta(\beta-1)^{-1}.
$$
Then we have
$$
\mu(S\G(\rr))=\mu(\G(\rr))\sum_{s\in S/\Gamma}\det(s)^{-\beta}
=\mu(\G(\rr))\zeta(\beta)\zeta(\beta-1)=1.\notag
$$
It follows that $S\G(\rr)$ is a subset of $\mtwo(\rr)$ of full
measure, and in view of $\ma=\glq\mr$, the subset
$\glq\G(\rr)=\G(\af)$ of $\ma$ has full measure. The scaling
condition $\mu(gY)=\det(g)^{-\beta}\mu(Y)$ completely  determines
$\mu$ on $\G(\af)$ by its restriction to $\G(\rr)$. Conversely, we
know from \cite[Lemma~2.4]{LLN} that any $\slr$-invariant measure on
$\G(\rr)$ extends uniquely to a measure on $\G(\af)$ satisfying the
scaling condition. To get a measure whose value on
$\mtwo(\rr)\cap\G(\af)$ is $1$, as required in Lemma~\ref{KMS1}, we
need the total mass of $\G(\rr)$ to be
$\zeta(\beta)^{-1}\zeta(\beta-1)^{-1}$. To summarize, the map
$$
\mu\mapsto \zeta(\beta)\zeta(\beta-1)\mu|_{\G(\rr)}
$$
defines a bijection between measures on $\ma$ satisfying the
conditions in Lemma~\ref{KMS1} and $\slr$-invariant probability
measures on $\G(\rr)$. In particular, extremal measures $\mu$
correspond to measures supported on one $\slr$-orbit.

Now fix $w\in\G(\rr)$ and consider the state $\varphi_{\beta,w}$
defined in Lemma~\ref{eextr}. Then
$$
\varphi_{\beta,w}(f)=\zeta(\beta)^{-1}\zeta(\beta-1)^{-1}
\sum_{s\in\Gamma\backslash S}\det(s)^{-\beta} f(sw)
$$
for $f\in C(\slr\backslash\mtwo(\rr))$. Notice that the sets $\slr
s$ are closed and disjoint for $s$ lying in different right cosets
of $\Gamma$ in~$S$, and their union is $\mtwo(\rr)\cap\glap$. We
thus see that the $\slr$-invariant measure~$\mu_{\beta,w}$ on
$\mtwo(\rr)$ defined by $\varphi_{\beta,w}$ has the property
$$
\mu_{\beta,w}(\slr sw)=\zeta(\beta)^{-1}\zeta(\beta-1)^{-1}
\det(s)^{-\beta} \text{ for } s\in S.
$$
Since $\sum_{s\in\Gamma\backslash S}
\zeta(\beta)^{-1}\zeta(\beta-1)^{-1} \det(s)^{-\beta}=1$, we
conclude that $(\mtwo(\rr)\cap\glap)w$ is a subset of $\mtwo(\rr)$
of full measure, so that $\varphi_{\beta,w}$ is regular at every
prime. Moreover, the restriction of $\mu_{\beta,w}$ to $\G(\rr)$ is
supported on one orbit $\slr w$. So indeed $\varphi_{\beta,w}$ is a
regular extremal KMS-state corresponding to the orbit $\slr w$, and
the proof of the theorem is complete. \ep

\begin{remark}\label{fin} \mbox{\ }
\enu{i} As we observed in the proof of Theorem~\ref{main_theorem},
for $\beta>1$ and every prime $p$ the state
$\varphi_{\beta,p}=\psi_{\beta,p}\circ\pi_p$ is the unique
$\sigma$-KMS$_\beta$-state on
$\mathfrak{A}_p^{\Z_p^*}={\redheck{\mtwo(\Z[p^{-1}]) \rtimes
G_p}{\mz\rtimes\Gamma}}^{\Z_p^*}$ which is normal with respect
to~$\pi_p$. Since $\varphi_{\beta,p}$ is a factor state, the state
$\otimes_{p\in\primes}\varphi_{\beta,p}$ on $\bigotimes_{p\in
\primes} \mathfrak{A}_p^{\Z_p^*}$ is the unique KMS$_\beta$-state
whose restriction to the factor corresponding to a prime $p$
coincides with $\varphi_{\beta,p}$. Thus by Theorem~\ref{Symm},
$\redheck{P}{P_0}^{\rr^*}$ has a unique $\sigma$-KMS$_\beta$-state
regular at every prime. So at least for regular KMS-states the
situation is similar to the one dimensional case: the group $\rr^*$
acts on the algebra, the fixed point algebra is a tensor product of
algebras corresponding to different primes, and for $\beta>1$ this
fixed point algebra has a unique regular $\sigma$-KMS$_\beta$-state,
which is a product-state. \enu{ii} We claim that there are no
$\sigma$-KMS$_\beta$-states on $\redheck{P}{P_0}$ for $\beta<0$.
Indeed, if $\varphi$ is such a state then for the isometry
$u_p=\left[\diag{p}{p}\right]\in\hecke{S}{\Gamma}$ we have
$$
1=\varphi(u^*_pu_p)=p^{2\beta}\varphi(u_pu^*_p)\le p^{2\beta},
$$
which is a contradiction. On the other hand, for $\beta>0$,
$\beta\ne1$, we conjecture that any $\sigma$-KMS$_\beta$-state on
$\redheck{P}{P_0}$ is regular at every prime. Indeed, for the full
Connes-Marcolli $\G$-system the analogous property for KMS-states
holds by \cite[Corollary~3.6]{LLN}. To prove regularity for states
on $\redheck{P}{P_0}$ it would be
enough to show that $\mathfrak{A}_p^{\Z^*_p}$ has no
$\sigma$-KMS$_\beta$-states for $\beta\in(0,1)$, and a unique
KMS$_\beta$-state for $\beta>1$, namely, $\varphi_{\beta,p}$. As we
remarked in the proof of Theorem~\ref{main_theorem}, to check
whether a KMS-state coincides with $\varphi_{\beta,p}$, it is enough
to evaluate it at one nonzero positive element $a_p$ such that
$\pi_p(a_p)$ is compact. E.g. we can take
$$
a_p=v^*_pv_p-v_pv_p^*-p(1-u_pu_p^*),
$$
which by \proref{projection} is the preimage of the rank one
projection onto $\C\delta_{\Z^2}\subset\ell^2(\Gamma\backslash
S_p)$. From this, one can conclude that $\varphi_{\beta,p}$ is
characterized by the equality
$$
\varphi_{\beta,p}(v_p^*v_p) 
=p+1.
$$
What makes the situation more difficult than the one dimensional
case is that the Toeplitz-Hecke algebra at prime $p$, that is, the
C$^*$-algebra generated by $\pi_p(\hecke{S_p}{\Gamma})$, does have
other KMS-states, which can be obtained, for example, by
representing the Toeplitz-Hecke algebra  on $\ell^2(\Gamma\backslash
S_p/\Gamma)$. So considerations in the Toeplitz-Hecke algebra alone
are not enough to prove the conjecture, and a better understanding of
the whole algebra $\mathfrak{A}_p^{\Z^*_p}$ is required.
\end{remark}

\end{document}